\documentclass[11pt]{amsart}

\usepackage{verbatim}
\usepackage{amsmath}
\usepackage{amssymb}
\usepackage{amsfonts}
\usepackage{amsthm}
\usepackage{mathrsfs}
\usepackage{tikz}
\usetikzlibrary{matrix,arrows}
\usepackage[all]{xy}
\usepackage{mathrsfs}
\usepackage{a4wide}
\usepackage{enumerate}
\usepackage{hyperref}

\theoremstyle{plain}
\newtheorem{theorem}[]{Theorem}
\newtheorem{lemma}[theorem]{Lemma}
\newtheorem{proposition}[theorem]{Proposition}
\newtheorem{corollary}[theorem]{Corollary}
\newtheorem{claim}[theorem]{Claim}

\newtheorem{conjecture}[theorem]{Conjecture}
\newtheorem*{thm}{Theorem}

\theoremstyle{definition}
\newtheorem{definition}[theorem]{Definition}
\newtheorem{example}[theorem]{Example}

\theoremstyle{remark}
\newtheorem{remark}[theorem]{Remark}

\newcommand{\G}{\mathcal{G}}
\newcommand{\E}{\widehat{E}}

\newcommand{\ucong}{\rotatebox{90}{$\cong$}}

\title{Spaces of Graphs, Boundary Groupoids and the Coarse Baum-Connes Conjecture} 
\date{July 2013}
\author{Martin Finn-Sell \and Nick Wright}
\email{Martin Finn-Sell: M.Finn-Sell@soton.ac.uk \\ Nick Wright: N.J.Wright@soton.ac.uk}

\begin{document}
\bibliographystyle{alpha}
\begin{abstract}
We introduce a new variant of the coarse Baum-Connes conjecture designed to tackle coarsely disconnected metric spaces called the boundary coarse Baum-Connes conjecture. We prove this conjecture for many coarsely disconnected spaces that are known to be counterexamples to the coarse Baum-Connes conjecture. In particular, we give a geometric proof of this conjecture for spaces of graphs that have large girth and bounded vertex degree. We then connect the boundary conjecture to the coarse Baum-Connes conjecture using homological methods, which allows us to exhibit all the current uniformly discrete counterexamples to the coarse Baum-Connes conjecture in an elementary way.
\end{abstract}

\maketitle

\section{Introduction and Outline}
The coarse Baum-Connes conjecture for metric spaces plays a central role in answering positively certain topological and group theoretic problems \cite{MR1388312}; it implies a positive solution to the Novikov conjecture for finitely generated discrete groups \cite{MR1779613,MR1388300} and information about the existence of positive scalar curvature metrics on spin manifolds \cite{MR1388300,MR1817560}. It is well known that this conjecture has counterexamples \cite{MR1817560,higsonpreprint,explg1} in the class of coarsely disconnected spaces.

In this paper we introduce a new conjecture associated to uniformly discrete bounded geometry metric spaces that is tailored to studying the geometry of coarsely disconnected spaces by utilising their asymptotic geometry. For a space $X$ this conjecture is phrased in terms of a certain Baum-Connes conjecture for a naturally constructed groupoid called the \textit{boundary groupoid} of $X$, which we denote by $G(X)|_{\partial\beta X}$. 
\begin{conjecture} [Boundary Coarse Baum-Connes Conjecture]
Let $X$ be a uniformly discrete bounded geometry metric space. Then the assembly map:
\begin{equation*}
\mu_{bdry}:K_{*}^{top}(G(X)|_{\partial\beta X}, A_{\partial}) \rightarrow K_{*}(A_{\partial}\rtimes_{r}G(X)|_{\partial\beta X})
\end{equation*}
is an isomorphism.
\end{conjecture}

The coefficient algebra $A_{\partial}$ in the conjecture above is related to the asymptotic geometry of $X$ and occurs in the work of Yu and Oyono-Oyono \cite{MR2568691}. The idea here is that this conjecture should be easier to prove in those instances that the coarse Baum-Connes conjecture fails, an idea that we confirm during the first half of this paper. 

In particular we prove the following result:

\begin{theorem}
The boundary coarse Baum-Connes conjecture holds for the following classes of coarsely disconnected spaces:
\begin{enumerate}
\item box spaces of finitely generated a-T-menable groups;
\item coarse disjoint unions constructed from sequences of finite graphs with large girth and uniformly bounded above regularity.
\end{enumerate}
\end{theorem}

The first point is an easy observation outlined in Section \ref{Sect:CE}, however the second point is rather more interesting and relies on carefully studying the natural geometry of large girth sequences to produce a \textit{partial action} of free group \cite{MR2041539}. These ideas are covered in Section \ref{Sect:S3}.

The second objective of this paper is to connect the boundary coarse Baum-Connes conjecture, via homological methods, to the coarse Baum-Connes conjecture. Using this machinery we give elementary proofs of many of the counterexample arguments from \cite{MR1911663,higsonpreprint,explg1} as well as many results concerning classes of expander graphs present in the literature \cite{MR2419930,MR2764895,MR2568691}; in particular in Section \ref{Sect:apps} we give an elementary proof of results of Willett and Yu \cite{explg1} concerning large girth sequences:

\begin{theorem}
Let $X$ be a space of graphs with large girth constructed from a sequence of finite graphs with uniformly bounded vertex degree. Then:
\begin{enumerate}
\item the coarse Baum-Connes assembly map is injective for $X$.
\item If $X$ comes from an expanding sequence then the coarse Baum-Connes assembly map for $X$ is not surjective.
\end{enumerate}
\end{theorem} 

Finally in Section \ref{Sect:Counter} we give a counterexample to the boundary coarse Baum-Connes conjecture by considering a space introduced by Wang \cite{MR2363697} and adaptations of the counterexample arguments present in the literature:

\begin{theorem}
There is a coarsely disconnected space $Y$ constructed from a box space of $SL_{2}(\mathbb{Z})$ for which the boundary coarse Baum-Connes assembly map is injective but fails to be surjective.
\end{theorem}

\subsection{Acknowledgements}
The first author wishes to thank the second author, his supervisor, for his unending patience and support as well as Rufus Willett for carefully reading the first drafts.
\numberwithin{theorem}{section}
\setcounter{theorem}{0}
\subsection{Groupoids, expanders and the coarse Baum-Connes conjecture.}\label{Sect:GO}
In this section we recall some definitions and results from the literature that will appear within the text. We begin with groupoids by focusing on certain relevant examples.

A groupoid is \textit{principal} if $(r,s): \G \rightarrow \G^{(0)} \times \G^{(0)}$ is injective and \textit{transitive} if $(r,s)$ is surjective. A groupoid $\G$ is a \textit{topological groupoid} if both $\G$ and $\G^{(0)}$ are topological spaces, and the maps $r,s, ^{-1}$ and the composition are all continuous. A Hausdorff, locally compact topological groupoid $\G$ is \textit{proper} if $(r,s)$ is a proper map and \textit{\'etale} or \textit{r-discrete} if the map $r$ is a local homeomorphism. When $\G$ is \'etale, $s$ and the product are also local homeomorphisms, and $\G^{(0)}$ is an open subset of $\G$.

\begin{definition}
Let $\G$ be a groupoid and let $x,y \in \G^{(0)}$ and $A,B \subset \G^{(0)}$. Set:
\begin{enumerate}
\item $\G_{x}=s^{-1}(x)$
\item $\G^{y}=r^{-1}(y)$
\item $\G^{y}_{x}=\G^{y} \cap \G_{x}$
\end{enumerate}
Denote by $\G|_{A}$ the subgroupoid $\G_{A}^{A}$, called the \textit{reduction} of $\G$ to $A$.
\end{definition}

\begin{definition}
Let $\G$ be a locally compact groupoid and let $Z$ be a locally compact space. $\G$ acts on $Z$ (or $Z$ is a $\G$-space) if there is a continuous, open map $r_{Z}: Z \rightarrow \G^{(0)}$ and a continuous map $(\gamma, z) \mapsto \gamma .z$ from $\G \ast Z:= \lbrace (\gamma, z) \in \G \times Z | s_{\G}(\gamma)=r_{Z}(z)\rbrace$ to $Z$ such that $r_{Z}(z).z=z$ for all $z$ and $(\eta \gamma).z= \eta.(\gamma. z)$ for all $\gamma, \eta \in \G^{(2)}$ with $s_{\G}(\gamma)=r_{Z}(z)$.
\end{definition}

When it is clear we drop the subscripts on each map. Right actions are dealt with similarly, replacing each incidence of $r_{Z}$ with $s_{Z}$.

\begin{definition}
Let $\G$ act on $Z$. The action is said to be \textit{free} if $\gamma.z=z$ implies that $\gamma = r_{Z}(z)$.
\end{definition}
We end this section with some useful examples.

\begin{example}\label{Ex:TransGrp}
Let $X$ be a topological $\Gamma$-space. Then the \textit{transformation groupoid} associated to this action is given by the data $X \times G \rightrightarrows X$ with $s(x,g)=x$ and $r(x,g)=g.x$. We denote this by $X \rtimes G$. A basis $\lbrace U_{i} \rbrace$ for the topology of $X$ lifts to a basis for the topology of $X \rtimes G$, given by sets $[U_{i},g]:=\lbrace (u,g) | u \in U_{i} \rbrace$. 
\end{example}

\begin{example}
The construction in the example above can be generalized to actions of \'etale groupoids. We are concerned with the topology here: Given an \'etale groupoid $\G$ and a $\G$-space $X$ as well as a with a basis $\lbrace U_{i} \rbrace$ for $\G^{(0)}$. We can pull this basis back to a basis for $X \rtimes \G$ given by $[r_{z}^{-1}(U_{i}),\gamma]$, where $U_{i} \subseteq s(\gamma)$.
\end{example}

\begin{example}Let $X$ be a coarse space with uniformly locally finite, weakly connected coarse structure $\mathcal{E}$. Define $G(X):=\cup_{E\in \mathcal{E}}\overline{E} \subseteq \beta(X\times X)$. Then $G(X)$ is a locally compact, Hausdorff topological space. To equip it with a product and inverse we would ideally have liked to be considering the natural extension of the pair groupoid product on $\beta X \times \beta X$.  We remark that the map $(r,s)$ from $X \times X$ extends first to an inclusion into $\beta X \times \beta X$ and universally to $\beta( X \times X)$, giving a map $(r,s): \beta (X\times X) \rightarrow \beta X \times \beta X$. We can restrict this map to each entourage $E \in \mathcal{E}$ allowing us to map the set $G(X)$ to $\beta X \times \beta X$. 
\end{example}

\begin{lemma}\label{Lem:CorRoe}\cite[Corollary 10.18]{MR2007488}
Let $X$ be a uniformly discrete bounded geometry metric space and let $E$ be any entourage. Then the inclusion $E \rightarrow X \times X$ extends to an injective homeomorphism $\overline{E} \rightarrow \beta X \times \beta X$, where $\overline{E}$ denotes the closure of $E$ in $\beta(X \times X)$.
\end{lemma}

Using this Lemma, we can conclude that the groupoid $G(X)$ can be embedded topologically into $\beta X \times \beta X$ and so we can equip it with the induced product and inverse. 

As we are considering the metric coarse structure we can reduce this to considering only generators:
\begin{equation*}
G(X):=\bigcup_{R>0}\overline{\Delta_{R}}
\end{equation*}

This groupoid plays an important role the groupoid formulation of the coarse Baum-Connes conjecture and the most general consequences of that conjecture. We recall some results concerning this groupoid from the literature:
\begin{theorem}
Let $X$ be a uniformly discrete bounded geometry metric space. Then following hold:
\begin{enumerate}
\item $G(X)$ is an \'etale locally compact Hausdorff principal topological groupoid with unit space $G(X)^{(0)}=\beta X$. \cite[Theorem 10.20]{MR2007488}\cite[Proposition 3.2]{MR1905840};
\item $C^{*}_{r}(G(X))$ is isomorphic to the uniform Roe algebra $C^{*}_{u}(X)$. \cite[Proposition 10.29]{MR2007488};
\item The coarse Baum-Connes conjecture for $X$ is equivalent to the Baum-Connes conjecture for $G(X)$ with coefficients in $\ell^{\infty}(X,\mathcal{K})$. \cite[Lemma 4.7]{MR1905840}.
\end{enumerate}
\end{theorem}

Lastly, we recall the definition of an expander:

\begin{definition}
Let $\lbrace X_{i} \rbrace$ be a sequence of finite graphs and let $X$ be the associated space of graphs. Then the space $X$ (or the sequence $\lbrace X_{i} \rbrace$) is an \textit{expander} if:
\begin{enumerate}
\item There exists $k\in \mathbb{N}$ such that all the vertices of each $X_{i}$ have degree at most $k$.
\item $\vert X_{i} \vert \rightarrow \infty$ as $i\rightarrow \infty$.
\item There exists $c>0$ such that $spectrum(\Delta_{i})\subseteq \lbrace 0 \rbrace \cup [c,1]$ for all $i$.
\end{enumerate}
\end{definition}

\begin{remark}\label{Rem:Ghost}
Each Laplacian $\Delta_{i}$ has propagation $1$, so we can form the product in the (algebraic) Roe algebra:
\begin{equation*}
\Delta:=\prod_{i}(\Delta_{i} \otimes q) \in C^{*}_{u}(X)\otimes \mathcal{K} \subset C^{*}X
\end{equation*}
Now we can consider projection $p$ onto the kernel of $\Delta$. For an expander $X$ we have $spectrum(\Delta) \subseteq \lbrace 0 \rbrace \cup [c,1]$ for some $c>0$, and so by an application of the functional calculus we can conclude that $p \in C^{*}X$. As $\Delta$ breaks up as a product we observe that its $Ker(\Delta)=\oplus_{i} Ker(\Delta_{i})$ and so the projection $p$ decomposes as a product:
\begin{equation*}
p=\prod_{i}p^{(i)}.
\end{equation*}
Additionally, it is easy to see that a function in $\ell^{2}(X_{i})$ is an element of the kernel of $\Delta_{i}$ if and only if it is a constant function, and so each $p^{i}$ has matrix entries $p^{i}_{x,y}=\frac{1}{\vert X_{i} \vert}$. 
\end{remark}
The following notion is due to Guoliang Yu (unpublished):

\begin{definition}\label{Def:Ghost}
An operator $T \in C^{*}X$ is a \textit{ghost operator} if $\forall \epsilon >0$ there exists a bounded subset $B \subset X\times X$ such that the norm: $\Vert T_{xy} \Vert \leq \epsilon$ for all $(x,y) \in (X\times X) \setminus B$.
\end{definition}

It is clear, from the definition as a product given above, that the kernel $p$ of the Laplacian $\Delta$ is a ghost operator.

\section{The Boundary Coarse Baum-Connes Conjecture.}\label{Sect:CE}

Throughout this section let $G$ be an \'etale, Hausdorff locally compact topological groupoid. For what we outline below these are essentially unnecessary assumptions, however our focus on the coarse groupoid leads us to study this class. We outline the main concept introduced in \cite{MR1911663}. 

\begin{definition}
A subset of $F\subseteq G^{(0)}$ is said to be \textit{saturated} if for every element of $\gamma \in G$ with $s(\gamma) \in F$ we have $r(\gamma) \in F$. For such a subset we can form subgroupoid of $G$, denoted by $G_{F}$ which has unit space $F$ and $G_{F}^{(2)}=\lbrace \gamma \in G | s(\gamma) \mbox{ and } r(\gamma) \in F \rbrace$.
\end{definition}

We will be considering closed saturated subsets $F$. Any subset of this type gives rise to an algebraic decomposition:
\begin{equation*}
G = G_{F^{c}}\sqcup G_{F}
\end{equation*}
This lets us construct maps on the $*-$algebras of compactly supported functions associated with $G$,$G_{F}$ and $G_{F^{c}}$:
\begin{equation*}
0 \rightarrow C_{c}(G_{F^{c}}) \rightarrow C_{c}(G) \rightarrow C_{c}(G_{F}) \rightarrow 0.
\end{equation*}
Where the quotient map $C_{c}(G) \rightarrow C_{c}(G_{F})$ is given by restriction and the inclusion $C_{c}(G_{F^{c}}) \rightarrow C_{c}(G)$ is given by extension. By the functorial properties of the maximal $C^{*}$-norm this extends to the maximal groupoid $C^{*}$-algebras:
\begin{equation*}
0 \rightarrow C_{max}^{*}(G_{F^{c}}) \rightarrow C_{max}^{*}(G) \rightarrow C_{max}^{*}(G_{F}) \rightarrow 0.
\end{equation*}
On the other hand this may fail to be an exact sequence when we complete in the norm that arises from the left regular representation $\lambda_{G}$. This is detected at the level of K-theory, as discussed in \cite{MR1911663}, by considering the sequence:
\begin{equation}\label{eqn:neim}
K_{0}(C^{*}_{r}(G_{U}))\rightarrow K_{0}(C^{*}_{r}(G)) \rightarrow K_{0}(C^{*}_{r}(G_{F}))
\end{equation}
This was used in \cite{MR1911663} to construct multiple different types of counterexample to the Baum-Connes conjecture for groupoids. We observe that whilst the sequence:
\begin{equation*}
\xymatrix{
0 \ar[r] & C_{r}^{*}(G_{F^{c}}) \ar[r]^{\alpha}& C_{r}^{*}(G) \ar[r]^{q} & C_{r}^{*}(G_{F}) \ar[r] & 0
}
\end{equation*}
may fail to be exact in the middle term the maps $\alpha$ and $q$ both exist and the map $q$ is surjective by considering the following diagram.
\begin{equation*}
\xymatrix{C^{*}_{max}(G) \ar@{->>}[r] \ar@{->>}[d] & C^{*}_{max}(G_{F}) \ar@{->>}[d]\\
C^{*}_{r}(G)\ar@{->}[r] &   C^{*}_{r}(G_{F})
}
\end{equation*}
It is also clear that the image of $\alpha$ is contained in the kernel of $q$, whence we can make the sequence exact artificially by replacing $C_{r}^{*}(G_{F^{c}})$ by the ideal $I:=\ker(q)$. We can then define a new assembly map in the first term to be the composition of the original assembly map $\mu_{F^{c}}$ and the K-theory map induced by inclusion $i_{*}:K_{*}(C^{*}_{r}(G_{F^{c}})) \rightarrow K_{*}(I)$. Then in terms of assembly maps this gives us a new commutative diagram:
\begin{equation*}\label{Fig:F1}
\xymatrix@=1em{
\ar[r] & K_{1}(C^{*}_{r}(G_{F}) \ar[r] & K_{0}(I) \ar[r]& K_{0}(C^{*}_{r}(G)) \ar[r]& K_{0}(C^{*}_{r}(G_{F})\ar[r] & K_{1}(I) \ar[r] & \\
\ar[r] & K_{1}^{top}(G_{F}) \ar[r] \ar[u]& K_{0}^{top}(G_{F^{c}}) \ar[r]\ar[u]& K_{0}^{top}(G) \ar[r]\ar[u]& K_{0}^{top}(G_{F}) \ar[r]\ar[u]& K_{1}^{top}(G_{F^{c}})\ar[u] \ar[r]& 
}
\end{equation*}
where the rows here are exact. 
As in \cite{MR1911663} we now consider $G=G(X)$, the coarse groupoid associated to some uniformly discrete bounded geometry metric space $X$.

\subsection{The Coarse Groupoid Conjecture}
Let $X$ be a uniformly discrete bounded geometry metric space. From what was described above we can associate to each closed saturated subset $F$ of the unit space space $\beta X$ a long exact sequence in K-theory. We consider the obvious closed saturated subset: $\partial \beta X \subset G(X)^{(0)}$. This gives us the following commutative diagram (omitting coefficients):
\begin{equation*}
\xymatrix@=0.7em{
 K_{1}(C^{*}_{r}(G(X)|_{\partial\beta X}) \ar[r] & K_{0}(I) \ar[r]& K_{0}(C^{*}_{r}(G(X))) \ar[r]& K_{0}(C^{*}_{r}(G(X)|_{\partial\beta X})\ar[r] & K_{1}(I)  \\
 K_{1}^{top}(G(X)|_{\partial\beta X}) \ar[r] \ar[u]& K_{0}^{top}(X \times X) \ar[r]\ar[u]& K_{0}^{top}(G(X)) \ar[r]\ar[u]& K_{0}^{top}(G(X)|_{\partial\beta X}) \ar[r]\ar[u]^{\mu_{bdry}}& K_{1}^{top}(X \times X)\ar[u]
}
\end{equation*}

We can now properly formulate the boundary conjecture by replacing the coefficients. Let $A_{\partial}:= \frac{\ell^{\infty}(X,\mathcal{K})}{C_{0}(X,\mathcal{K})}$

\begin{conjecture} [Boundary Coarse Baum-Connes Conjecture]
Let $X$ be a uniformly discrete bounded geometry metric space. Then the assembly map:
\begin{equation*}
\mu_{bdry}:K_{*}^{top}(G(X)|_{\partial\beta X}, A_{\partial}) \rightarrow K_{*}(A_{\partial}\rtimes_{r}G(X)|_{\partial\beta X})
\end{equation*}
is an isomorphism.
\end{conjecture}

\subsection{Tools to Prove the Boundary Conjecture.}
We consider the situation when a group $\Gamma$ acts on a space of graphs $X$. We introduce a definition:

\begin{definition}
Let $\mathcal{E}$ be a coarse structure on $X$ and let $\mathcal{S}$ be a family of subsets $\mathcal{E}$. We say that $S$ generates $\mathcal{E}$ at infinity if for all $E \in \mathcal{E}$: 
\begin{equation*}
E \subseteq (\bigcup_{k=1}^{n}S_{k})\cup F
\end{equation*}
Where each $S_{k} \in \mathcal{S}$ and $F$ is a finite subset of $X \times X$.
\end{definition}

\begin{remark}
The above definition is equivalent to asking that $\overline{E}\setminus E \subseteq \bigcup_{k=1}^{n}\overline{S_{k}}\setminus S_{k}$, where the closure takes place in $\beta X \times \beta X$.
\end{remark}

Recall that $\Delta_{g}:= \lbrace (x,x.g)|x \in X \rbrace$ is the $g$-diagonal in $X$.

\begin{proposition}\label{Prop:Crit}
Let $X$ be uniformly discrete bounded geometry metric space and let $\Gamma$ be a finitely generated discrete group. If $\Gamma$ acts on $X$ such that the induced action on $\beta X$ is free on $\partial \beta X$ and the action generates the metric coarse structure at infinity. Then $G(X)|_{\partial \beta X} \cong \partial \beta X \rtimes \Gamma$.
\end{proposition}
\begin{proof}
Recall that under the hypothesis that the $\Gamma$ action generates the metric coarse structure at infinity we know that $G(X) = (\bigcup_{g} \overline{\Delta_{g}})\cup (X \times X)$.

We consider a map from transformation groupoid $\beta X \rtimes \Gamma$ to $G(X)$. Observe that $\Delta_{g}$ is the bijective image of the set $\lbrace (x,g)|x \in X \rbrace$. We extend this map to the respective closures, giving a map for each $g$ from the set of $\lbrace (\omega,g) | \omega \in \beta X \rbrace$ to $\overline{\Delta_{g}}$. Let $\lbrace x_{\lambda} \rbrace$ be a net of elements in $X$ that converge to $\omega$. Then clearly each pair $(x_{\lambda},g)$ is mapped to $(x_{\lambda},x_{\lambda}.g)$ under the identification. Consider the element of the closure $\gamma_{g}= \lim_{\lambda}(x_{\lambda},x_{\lambda}.g) \in \overline \Delta_{g}$ and define the extension of the bijection to be the map: $(\omega,g) \mapsto \gamma_{g}$.

This map is well-defined as $G(X)$ is principal and so $\gamma_{g}$ is completely determined by its source and range - in particular any other net $y_{\lambda}$ that converges to $\omega$ gives rise to the same element $\gamma_{g}$. We can then extend this over the entire groupoid $\beta X \rtimes \Gamma$ element wise, where it is certainly continuous but not in general injective or surjective (injectivity would require a free action and surjectivity a transitive one). 

The map is a groupoid homomorphism because $G(X)$ is principal; this follows as principality implies the following diagram commutes:
\begin{equation*}
\xymatrix{
& \beta X \rtimes \Gamma \ar@/_/[rd]_{(r,s)} \ar[r] & G(X) \ar@{^{(}->}[d]^{(r,s)}\\ & & \beta X \times \beta X
}
\end{equation*}

We now restrict this map to the boundary $\partial\beta X$. We prove that the coronas $\overline{\Delta_{g}}\setminus \Delta_{g}$ are disjoint: let $\gamma \in \overline{\Delta_{g}}\setminus \Delta_{g}$ and $\overline{\Delta_{h}}\setminus \Delta_{h}$ for $g \not = h$. Then we have that $s(\gamma)=\omega, r(\gamma)=\omega.g=\omega.h$, hence we have a fixed point on the boundary, which is a contradiction. Consider the restricted diagram:
\begin{equation*}
\xymatrix{
& \partial\beta X \rtimes \Gamma \ar@{^{(}->}[rd]_{(r,s)} \ar[r] & G(X)_{\partial\beta X} \ar@{^{(}->}[d]^{(r,s)}\\ & & \partial\beta X \times \partial\beta X.
}
\end{equation*}
As the coronas are disjoint, for any $\gamma \in G(X)$ it is possible to find a unique $g$ and $\omega$ such that $s(\gamma)=\omega, r(\gamma)=\omega.g$. It then follows that the pair $(\omega,g)$ map onto $\gamma$. Hence the map is a surjection. To see injectivity, we appeal again to freeness of the action. The action being free implies that the groupoid $\partial\beta X \rtimes \Gamma$ is principal. From a brief consideration of the diagram above injectivity follows.

It remains to consider the topology. For each $g\in \Gamma$, the map from $\partial \beta X \rtimes \Gamma$ to $G(X)_{\partial\beta X}$ on the piece $[\partial\beta X,g]$ is a homeomorphism as $\overline{\Delta_{g}}\setminus\Delta_{g}$ can be identified with $\partial\beta X$, induced by the source map or range map. Hence the map takes clopen sets to clopen sets; these form a basis for the topology of $G(X)|_{\partial\beta X}$ and so the map is a global homeomorphism.
\end{proof}

\begin{example}(Box spaces)
Let $\Gamma$ be a finitely generated residually finite group and let $\lbrace N_{i}\rbrace$ be a family of finite index normal subgroups such that $N_{i}\leq N_{i+1}$, $\bigcap_{i \in \mathbb{N}}N_{i}=1$ and a fixed generating set $S$. Then the sequence of groups $\lbrace \Gamma/N_{i}\rbrace$ with generating sets $\pi_{i}(S)$ admits an action via quotient maps. Let $\square\Gamma= \sqcup_{i \in \mathbb{N}}\frac{\Gamma}{N_{i}}$, equipped with a metric that restricts to the metric induced from the generating sets $\pi_{i}(S)$ for each $i$, and has the property that $d(\frac{\Gamma}{N_{i}},\frac{\Gamma}{N_{j}})\rightarrow \infty$ as $i+j \rightarrow \infty$. This is called a \textit{box space} for $\Gamma$. The Stone-\u{C}ech boundary admits a free action of the group and the metric structure is generated at infinity by the quotients maps and the right action of the group via these maps. Proposition \ref{Prop:Crit} then provides a description of the boundary groupoid and converts the boundary conjecture into a case of the Baum-Connes conjecture with coefficients for $\Gamma$.

This process did not require a normal subgroup; box spaces can constructed using Schreier quotients. The conditions on the subgroups change to reflect the absence of normality. Let $\Gamma$ be a residually finite group and let $\lbrace H_{i} \rbrace$ be a family of nested subgroups of finite index with trivial intersection, and additionally satisfying: each $g \in \Gamma$ belongs to only finitely many conjugates of the subgroups from the family $\lbrace H_{i} \rbrace$. Fixing a left invariant metric on $\Gamma$ the box space can be constructed using the left quotients of $\Gamma$ by the $H_{i}$. In this instance these spaces are graphs with no left action of $\Gamma$. However they do retain a right action of $\Gamma$ that determines the metric at infinity and becomes free on the boundary (this is due to the additional constraint concerning conjugates of the $H_{i}$).
\end{example}

\section{A Proof of the Boundary Conjecture for Spaces of Graphs with Large Girth.}\label{Sect:MR}
The aim of this section is to prove the following result:
\begin{thm}
The boundary coarse Baum-Connes conjecture holds for spaces of graphs with large girth and uniformly bounded vertex degree.
\end{thm}

For this result we will need a weaker notion of action; We construct a partial action of $F_{k}$ for some finite $k$. Using a stronger version of Proposition \ref{Prop:Crit} we will show that the boundary groupoid $G(X)|_{\partial\beta X}$ is a transformation groupoid $\partial\beta X \rtimes \G_{\widehat{X}}$ for some \'etale groupoid $\G_{\widehat{X}}$. As before, this converts the boundary Baum-Connes conjecture into a special case of Baum-Connes with coefficients for the groupoid $\G_{\widehat{X}}$. Lastly, we construct a continuous, proper groupoid homomorphism from $\G_{\widehat{X}}$ to $F_{k}$, which we use to transfer the Haagerup property from $F_{k}$ to the groupoid $\G_{\widehat{X}}$. This will complete the proof.

As noted above we want to consider partial actions of a group $\Gamma$ on a space $X$. This means that the elements of $\Gamma$ give rise to partial bijections of $X$, i.e bijections between subsets of $X$. These partial bijections are algebraically encoded within an inverse semigroup, and the partial action of $\Gamma$ is a \textit{dual prehomomorphism} from $\Gamma$ into that semigroup.

\subsection{An Interlude into Inverse Semigroup Theory}\label{Sect:SemiToGpoid}
                                                                               
\begin{definition}\label{Def:invsemi}
Let $S$ be a semigroup. We say $S$ is $inverse$ if there exists a unary operation $*:S \rightarrow S$ satisfying the following identities:
\begin{enumerate}
\item $(s^{*})^{*}=s$
\item $ss^{*}s=s$ and $s^{*}ss^{*}=s^{*}$ for all $s \in S$
\item $ef=fe$ for all idempotents $e,f \in S$ 
\end{enumerate}
\end{definition}

Recall that a semigroup with a unit element is called a \textit{monoid}, and it makes sense to talk about inverse monoids in the obvious way. A very fundamental example is the \textit{symmetric inverse monoid} on any set $X$; consider the collection of all partial bijections of $X$ to itself, giving them them the natural composition law associated to functions, finding the largest possible domain.

When $X$ is a metric space we will be considering a inverse submonoid of $I(X)$ in which every partial bijection that maps elements only a finite distance, that is a generalised (or partial) translation. We denote this by $I_{b}(X)$.

\begin{definition}
Let $S$ be an inverse monoid. We denote by $E(S)$ the semilattice of idempotents (just by $E$ if the context is clear). This is a meet semilattice, where the meet is given by the product of $S$ restricted to $E$. In this situation, we can use the following partial order:
\begin{equation*}
e \leq f \Leftrightarrow ef=e
\end{equation*}
This order can be extended naturally to the entire of $S$: $s \leq t$ if there exists $e \in E(S)$ such that $s=et$. In terms of partial bijections this order corresponds to restricting an element to a subset of its domain.
\end{definition}

We remark that for a metric space $X$ every idempotent element in $I(X)$ moves elements no distance, and hence $E(I(X))=E(I_{b}(X))$.

We want to consider quotient structures of an inverse monoid. In general quotients are given by equivalence relations that preserve the semigroup structure.

\begin{definition}
An equivalence relation $\sim$ on $S$ is called a \textit{congruence} if for every $u,v,s,t \in S$ such that $s \sim t$, we know that $su\sim tu$ and $vs \sim vt$. This allows us to equip the quotient $\frac{S}{\sim}$ with a product, making it into an inverse monoid.
\end{definition}

One such example of this arises from an \textit{ideal} in $S$.

\begin{definition}
Let $I$ be a subset of $S$. $I$ is an ideal of $S$ if $SI \cup IS \subset I$.
\end{definition}

From an ideal we can get a quotient - at the cost of a \textit{zero element}, that is an element $0$ such that $0s=s0=0$ for all $s \in S$.

\begin{definition}
Let $S$ be an inverse monoid and let $I$ be an ideal of $S$. Then we can define $\frac{S}{I}$ to be the quotient of the set $S$ by the congruence: $x \sim y$ if $x=y$ or $x$ and $y$ are elements of $I$. 
\end{definition}

Another specific congruence we will be interested in is called the minimum group congruence on $S$.  This congruence, denoted by $\sigma$, is given by:
\begin{equation*}
s \sigma t \Leftrightarrow (\exists e \in E) es = et
\end{equation*}
This congruence is \textit{idempotent pure}, that is for $e \in E(S)$ and $s \in S$, $e \sim s$ implies $s \in E$.

\begin{definition}
We say $S$ is 0-E-unitary if $\forall e \in E\setminus 0, s \in S \setminus 0$ $e \leq s$ implies $s \in E$. We say it is 0-F-inverse if in addition there exists a subset $T \subset S$ such that for every $s \in S$ there exists a unique $t \in T$ such that $s \leq t$ and if for all $u \in S$ with $s \leq u$ then $u \leq t$.
\end{definition}

\subsection{Groupoids from Inverse Monoids}
We take an inverse monoid $S$ and produce a universal groupoid $\G_{\E}$. One way to do this involves studying the actions of $S$ on its semilattice $E$.

We outline the steps in the construction.
\begin{enumerate}
\item Build an action of $S$ on $E$.
\item Build a dual space to $E$, which is compact and Hausdorff. This is a \textit{Stone dual} to $E$. Show this admits an action of $S$.
\item Build the groupoid $\G_{\E}$ from this data.
\end{enumerate}

\begin{definition}
\begin{enumerate}
\item Let $D_{e}=\lbrace f \in E | f \leq e \rbrace$. For $ss^{*} \in E$, we can define a map $\rho_{s}(ss^{*})=s^{*}s$, extending to $D_{ss^{*}}$ by $\rho_{s}(e) = s^{*}es$. This defines a partial bijection on $E$ from $D_{ss^{*}}$ to $D_{s^{*}s}$. 

\item We consider a subspace of $\textbf{2}^{E}$ given by the functions $\phi$ such that $\phi(0)=0$ and $\phi(ef)=\phi(e)\phi(f)$. This step is a generalisation of Stone duality \cite{MR2672179}. We can topologise this as a subspace of $\textbf{2}^{E}$, where it is closed. This makes it compact Hausdorff, with a base of topology given by $\widehat{D}_{e}= \lbrace \phi \in \E | \phi(e)=1 \rbrace$. This admits a dual action induced from the action of $S$ on $E$. This is given by the pointwise equation for every $\phi \in \widehat{D}_{s^{*}s}$:
\begin{equation*}
\widehat{\rho}_{s}(\phi)(e)=\phi(\rho_{s}(e))=\phi(s^{*}es)
\end{equation*}
The use of $\widehat{D}_{e}$ to denote these sets is not a coincidence, as we have the following map $D_{e} \rightarrow \widehat{D}_{e}$:
\begin{equation*}
e \mapsto \phi_{e}, \phi_{e}(f)=1 \mbox{ if } e \leq f \mbox{ and } 0 \mbox{ otherwise }.
\end{equation*}
\begin{remark}
These character maps $\phi: E \rightarrow \lbrace 0,1 \rbrace$ have an alternative interpretation, they can be considered as \textit{filters} on $E$. A filter on $E$ is given by a set $F \subset E$ with the following properties:
\begin{itemize}
\item for all $e,f \in F$ we have that $e\wedge f=ef \in F$
\item for $e\in F$ with $e \leq f$ we have that $f \in F$ and
\item $0 \not\in F$
\end{itemize}
the relationship between characters and filters can be summarised as: To each character $\psi$ there is a filter:
\begin{equation*}
F_{\psi}= \lbrace e \in E | \psi(e)=1 \rbrace.
\end{equation*}
And every filter $F$ provides a character by considering $\chi_{F}$, its characteristic function.
\end{remark}

\item We take the set $S \times \E$, topologise it as a product and consider subset $\Omega:= \lbrace (s, \phi) | \phi \in \widehat{D}_{s^{*}s} \rbrace$ in the subspace topology. We then quotient this space by the relation:
\begin{equation*}
(s, \phi) \sim (t, \phi^{'}) \Leftrightarrow \phi=\phi^{'} \mbox{ and } (\exists e \in E) \mbox{ with } \phi \in \widehat{D}_{e} \mbox{ such that } es=et
\end{equation*}
We can give the quotient $\G_{\E}$ a groupoid structure with the product set, unit space and range and source maps:
\begin{eqnarray*}
\G_{\E}^{(2)}:=\lbrace ([s,\phi],[t,\phi^{'}]) | \phi=\widehat{\rho}_{t}(\phi^{'}) \rbrace \\
\G_{\E}^{(0)}:= \lbrace [e,\phi] | e \in E \rbrace \cong \E \\
s([t,\phi])=[t^{*}t,\phi], r([t,\phi])=[tt^{*},\phi], 
\end{eqnarray*}
and product and inverse:
\begin{eqnarray*}
[s,\phi].[t,\phi^{'}]= [st,\phi^{'}] \mbox{ if } ([s,\phi],[t,\phi^{'}]) \in \G_{\E}^{(2)}, [s,\phi]^{-1} = [s^{*},\widehat{\rho}_{s}(\phi)] 
\end{eqnarray*}
For all the details of the above, we refer to \cite[Section 4]{MR2419901}. This is the \textit{universal groupoid} associated to $S$. We collect some information about this groupoid from \cite{MR2419901,MR1724106} in Theorem \ref{Thm:Info}.
\end{enumerate}
\end{definition}

\begin{theorem}\label{Thm:Info}
Let $S$ be a countable 0-E-unitary inverse monoid, $E$ its semilattice of idempotents and $\G_{\E}$ its universal groupoid. Then the following hold for $\G_{\E}$:
\begin{itemize}
\item $\E$ is a compact, Hausdorff and second countable space.
\item $\G_{\E}$ is a Hausdorff groupoid \cite[Corollary 10.9]{MR2419901}.
\item Every representation of $S$ on Hilbert space gives rise to a covariant representation of $\G_{\E}$ and vice versa \cite[Corollary 10.16]{MR2419901}.
\item We have $C^{*}_{r}(S) \cong C^{*}_{r}(\G_{\E})$ \cite{MR1724106,MR1900993}
\end{itemize}
\end{theorem}

We make use of the following technical property that arises from the presence of maximal elements:

\begin{claim}\label{Claim:C1}
Let $S$ be 0-F-inverse. Then every element $[s,\phi] \in \G_{\E}$ has a representative $[t,\phi]$ where $t$ is a maximal element.
\end{claim}
\begin{proof}
Take $t=t_{s}$ the unique maximal element above $s$. Then we know 
\begin{equation*}
s = t_{s}s^{*}s \mbox{ and } s^{*}s \leq t_{s}^{*}t_{s}
\end{equation*} 
The second equation tells us that $t_{s}^{*}t_{s} \in F_{\phi}$ as filters are upwardly closed, thus $(t_{s},\phi)$ is a valid element. Now to see $[t_{s},\phi]=[s,\phi]$ we need to find an $e \in E$ such that $e \in F_{\phi}$ and $se=t_{s}e$. Take $e=s^{*}s$ and then use the first equation to see that $s(s^{*}s)=t_{s}(s^{*}s)$.
\end{proof}

\subsection{(Dual) Prehomomorphisms and General Partial Actions}\label{Sect:S3}
\begin{definition}
Let $\rho: S \rightarrow T$ be a map between inverse semigroups. This map is called a prehomomorphism if for every $s,t \in S$, $\rho(st) \leq \rho(s)\rho(t)$ and a dual prehomomorphism if for every $s,t \in S$ $\rho(s)\rho(t) \leq \rho(st)$.
\end{definition}

We recall that a congruence is said to be \textit{idempotent pure} for any $e \in E(S)$, $s \in S$ we have that $s$ is related to $e$ implies that $s \in E(S)$. We extend this definition to general maps in the following way.
\begin{definition}
A (dual) prehomomorphism $\rho$ is called \textit{idempotent pure} if $\rho(s)^{2}=\rho(s)$ implies $s \in E$.  
\end{definition}
In addition we call a map $S \rightarrow T$ \textit{0-restricted} if the preimage of $0 \in T$ is $0 \in S$.

\begin{definition}
Let $S$ be a 0-E-unitary inverse monoid. We say $S$ is \textit{strongly 0-E-unitary} if there exists an idempotent pure, 0-restricted prehomomorphism, $\Phi$ to a group $G$ with a zero element adjoined, that is: $\Phi:S \rightarrow G^{0}$. We say it is \textit{strongly 0-F-inverse} if it is 0-F-inverse and strongly 0-E-unitary. This is equivalent to the fact that the preimage of each group element under $\Phi$ contains a maximum element.
\end{definition}

\begin{example}
In \cite{MR745358,MR2221438} the authors introduce an inverse monoid that is universal for dual prehomomorphisms from a general inverse semigroup. In the context of a group $G$ This is called the \textit{prefix expansion}; its elements are given by pairs: $(X,g)$ for $\lbrace 1,g\rbrace \subset X$, where $X$ is a finite subset of $G$. The set of such $(X,g)$ is then equipped with a product and inverse:
\begin{equation*}
(X,g)(Y,h) = (X\cup gY,gh)\mbox{ , } (X,g)^{-1}=(g^{-1}X,g^{-1})
\end{equation*}
This has maximal group homomorphic image $G$, and it has the universal property that it is the largest such inverse monoid. We denote this by $G^{Pr}$. The partial order on $G^{Pr}$ can be described by reverse inclusion, induced from reverse inclusion on finite subsets of $G$. It is F-inverse, with maximal elements: $\lbrace(\lbrace 1,g \rbrace, g):g \in G \rbrace$.
\end{example}

\begin{definition}
Let $G$ be a finitely generated discrete group and let $X$ be a (locally compact Hausdorff) topological space. A \textit{partial action} of $G$ on $X$ is a dual prehomomorphism $\theta$ of $G$ in the symmetric inverse monoid $\mathcal{I}(X)$ that has the following properties:
\begin{enumerate}
\item The domain $D_{\theta_{g}^{*}\theta_{g}}$ is an open set for every $g$.
\item $\theta_{g}$ is a continuous map.
\item The union: $\bigcup_{g \in G}D_{\theta_{g}^{*}\theta_{g}}$ is $X$.
\end{enumerate}
\end{definition}

Given this data we can generate an inverse monoid $S$ using the set of $\theta_{g}$. This would then give a representation of $S$ into $\mathcal{I}(X)$. If the space $X$ is a coarse space, then it makes sense to ask if each $\theta_{g}$ is a close to the identity. In this case, we would get a representation into the bounded symmetric inverse monoid $\mathcal{I}_{b}(X)$. We call such a $\theta$ a \textit{bounded partial action} of $G$.

Let $X= \sqcup X_{i}$ be a space of graphs admitting a bounded partial action of a discrete group $G$. We remark that in this setting partial bijections in the group can have the following form: 
\begin{equation*}
\theta_{g}=\theta_{g}^{0} \sqcup \bigsqcup_{i>i_{0}}\theta_{g}^{i}.
\end{equation*} 
Where $i_{0}$ is the first $i$ for which the distances between the $X_{i}$s is greater than the upper bound of the distance moved by $\theta_{g}$, and the $\theta_{g}^{i}$ are componentwise partial bijections of the $X_{i}$. We collect all the additional pieces that act only between the first $i_{0}$ terms into $\theta_{g}^{0}$, which could be the empty translation. We remark now that it is possible that there are partial bijections $\theta_{g}$ that could have finite support, that is only finitely many terms that are non-empty after $i_{0}$. To avoid this we observe the following:

\begin{proposition}
Let $S = \langle \theta_{g} | g \in G \rangle$ and let $I_{fin}= \lbrace \theta_{g} | supp(\theta_{g}) \mbox{ is finite} \rbrace$. Then $I_{fin}$ is an ideal and the Rees quotient $S_{inf}=\frac{S}{I}$ is an inverse monoid with $0$.
\end{proposition}
\begin{proof}
To be an ideal, it is enough to show that $I_{fin}S \subset I_{fin}, SI_{fin}\subset I_{fin}$. Using the description of the multiplication of partial bijections from section \ref{Sect:SemiToGpoid}, it is clear that either combination $si$ or $is$ yields an element of finite support. Now we can form the Rees quotient, getting an inverse monoid with a zero - the zero element being the equivalence class of elements with finite support.
\end{proof}

We want to utilise a partial action to construct a groupoid, so we apply the general construction outlined in section \ref{Sect:SemiToGpoid} to get an improved version of Proposition \ref{Prop:Crit}. To do this requires a notion of translation length:
\begin{definition}
The length of each $\theta_{g}$ is defined to be:
\begin{equation*}
\vert \theta_{g} \vert = \sup \lbrace d(x,\theta_{g}(x)) : x \in Dom(\theta_{g})\rbrace.
\end{equation*} 
\end{definition}

\begin{definition}
Recall that we say a bounded partial action  $\theta$ generates the metric coarse structure at infinity if for all $R>0$ there exists $S>0$ such that $\overline{\Delta_{R}}\setminus \Delta_{R} \subseteq \bigcup_{\vert \theta_{g} \vert < S}\overline{\Delta_{\theta_{g}}}\setminus \Delta_{\theta_{g}}$. We say it finitely generates the metric coarse structure if the number of $\theta_{g}$ required for each $R$ is finite.
\end{definition}

\begin{remark}
Recall a groupoid $G$ is said to be \textit{principal} if the map $(s,r): G \rightarrow G^{(0)}\times G^{(0)}$ is injective.
\end{remark}

\begin{proposition}\label{Prop:Aug}
Let $\lbrace X_{i} \rbrace$ be a sequence of finite graphs and let $X$ be the corresponding space of graphs. If $\theta:G \rightarrow \mathcal{I}(X)$ is a bounded partial action of $G$ on $X$ such that the induced action on $\beta X$ is free on $\partial \beta X$, the inverse monoid $S_{inf}$ is 0-F-inverse with maximal element set $\lbrace \theta_{g} |g \in G\rbrace$ and the partial action finitely generates the metric coarse structure at infinity then there is a second countable, \'etale topological groupoid $\G_{\widehat{X}}$ such that $G(X)|_{\partial\beta X} \cong \partial\beta X \rtimes \G_{\widehat{X}}$.
\end{proposition}
\begin{proof}
Observe now that the finite $\theta_{g}$ play no role in the action on the boundary and so we work with $S_{inf}$. We build the groupoid from the bottom up, by first constructing the unit space using Proposition 10.6 and Theorem 10.16 from \cite{MR2419901}.

We consider the representation of the inverse monoid $S_{inf}$ on $\ell^{2}(X)$ induced by $\theta$ to get a representation $\pi_{\theta}:S \rightarrow \mathcal{B}(\ell^{2}(X))$. We can complete the semigroup ring in this representation to get an algebra $C^{*}_{\pi_{\theta}}S$, which has a unital commutative subalgebra $C^{*}_{\pi_{\theta}}E$. Proposition 10.6 \cite{MR2419901} then tells us that the spectrum of this algebra, which we will denote by $\widehat{X}$, is a subspace of $\E$ that is closed and invariant under the action of $S$ on which the representation $\pi_{\theta}$ is supported.

As the space $\widehat{X}$ is closed and invariant we can reduce the universal groupoid $\G_{\E}$ for $S_{inf}$ to $\widehat{X}$. This we denote by $\G_{\widehat{X}}$. 

We show this groupoid acts on $\beta X$; we make use of the assumption that $\theta_{g}$ is a bounded partial bijection for each $g \in G$ and again of the representation $\pi_{\theta}$. Each $\theta_{g}$ bounded implies that the algebra $C^{*}_{\theta}(S)$ is a subalgebra of $C^{*}_{u}X$. We now remark that the representation $\pi_{X}$, when restricted to $C^{*}E$ assigns each idempotent a projection in $C^{*}_{u}X$, that is $C^{*}_{\pi_{X}}(E)=\pi_{X}(C^{*}E) \subset \ell^{\infty}(X)$. Taking the spectra associated to this inclusion then gives us a map:
\begin{equation*}
r_{\beta X}: \beta X \twoheadrightarrow \widehat{X}
\end{equation*}
which is continuous. In particular as both $\beta X$ and $\widehat{X}$ are compact Hausdorff spaces; this map is closed (and open) and hence a quotient. We make use of this to define an action on $\beta X$. By Claim \ref{Claim:C1} we have that each element of our groupoid $\G_{\widehat{X}}$ can be represented by a pair $[\theta_{g}, \phi]$, for some $\phi \in \widehat{X}$. Observe also that as $X$ is discrete so are all of it subspaces, hence the maps $\theta_{g}$ are continuous (open) for each $g \in G$. These then extend to $\beta X$, and so coupled with the map $r_{\beta X}$ we can act by:
\begin{equation*}
[\theta_{g},\phi].\omega = \theta_{g}(\omega)
\end{equation*}
for all $\omega \in D_{\theta_{g}^{*}\theta_{g}}$ with $r_{Z}(\omega) = \phi$. We see that $r_{Z}(\omega).\omega= [\theta_{e}, r_{Z}(\omega)].\omega = \omega.$ and for all $([\theta_{g},\phi],[\theta_{h},\phi^{'}]) \in \G_{\widehat{X}}^{(2)}$ with $\phi^{'}=r_{Z}(\omega)$ we have:
\begin{equation*}
[\theta_{g}\theta_{h},\phi^{'}].\omega = \theta_{g}\theta_{h}(\omega) = \theta_{g}([\theta_{h},\phi^{'}].\omega)=[\theta_{g},\phi].([\theta_{h},\phi^{'}].\omega)
\end{equation*}
as $r_{Z}([\theta_{h},\phi^{'}].\omega)=\theta_{h}(\phi^{'})=\phi$.

It remains to prove the isomorphism of topological groupoids: $G(X)|_{\partial\beta X} \cong \partial\beta X \rtimes \G_{\widehat{X}}$. We follow the scheme of Proposition \ref{Prop:Crit} and build a map from $\beta X\rtimes \G_{\widehat{X}}$ to $G(X)$. Recall that as the partial action of $G$ generates the metric coarse structure at infinity $G(X)= (\bigcup_{g} \Delta_{\theta_{g}})\cup (X \times X)$. We observe that each $\Delta_{\theta_{g}}$ maps bijectively onto the domain of $\theta_{g}$, a subset of $X$.

This map extends to the closure of the domain precisely as in Proposition \ref{Prop:Crit}, where here we map the pair $(\omega,[\theta_{g},\phi])$ to the element $\gamma_{g,\phi}$ that is the limit $\lim_{\lambda}(x_{\lambda},\theta_{g}(x_{\lambda}))$ for some net $\lbrace x_{\lambda} \rbrace$ that converges to $\omega$ (and also to $\phi$). This map is well defined as the groupoid $G(X)$ is principal, and it fits into the following commutative diagram:
\begin{equation*}
\xymatrix{
& \beta X \rtimes \Gamma \ar@/_/[rd]_{(r,s)} \ar[r] & G(X) \ar@{^{(}->}[d]^{(r,s)}\\ & & \beta X \times \beta X
}
\end{equation*}
Again by principality, we can deduce that the covering map is a groupoid homomorphism.

We now restrict this map to the boundary $\partial\beta X$. As we know that the group action generates the metric coarse structure at infinity and that the partial action of the group $G$ is free on the boundary. Using these facts we can see that:
\begin{enumerate}
\item $\partial\beta X \rtimes \G_{\widehat{X}}$ is principal.
\item $G(X)|_{\partial\beta X} = \bigsqcup_{g}\overline{\Delta_{\theta_{g}}}\setminus \Delta_{\theta_{g}}$. 
\end{enumerate}
From both (1) and (2) we can further deduce that the covering map is a bijection on the boundary. Both groupoids are also \'etale and so each component $\overline{\Delta_{\theta_{g}}}\setminus \Delta_{\theta_{g}}$ is mapped homeomorphically onto its image and is therefore clopen. It follows then that we get the desired isomorphism $\partial\beta X \rtimes \G_{\widehat{X}} \cong G(X)|_{\partial\beta X}$ of topological groupoids.\end{proof}

We are interested in understanding the analytic properties of the groupoid $\G_{\widehat{X}}$. In particular we are interested in showing that the groupoid has the Haagerup property. To do this we study the inverse monoid $S$ associated to the partial action $\theta$.

\begin{proposition}\label{Prop:Strongly}
Let $S = \langle \theta_{g} | g \in G \rangle$, where $\theta: G \rightarrow S$ is a dual prehomomorphism. If $S$ is 0-F-inverse with $Max(S) = \lbrace \theta_{g} | g \in G \rbrace$ where each nonzero $\theta_{g}$ is not idempotent when $g \not = e$ then $S$ is strongly 0-F-inverse.
\end{proposition}
\begin{proof}
We build a map $\Phi$ back onto $G^{0}$. Let $m: S\setminus \lbrace 0 \rbrace \rightarrow Max(S)$ be the map that sends each non-zero $s$ to the maximal element $m(s)$ above $s$ and consider the following diagram:
\begin{equation*}
\xymatrix{
G\ar@{->}[r]^{\theta}\ar@{->}[dr]^{}  & S\ar@{->}[dr]^{\Phi}  & \\
  & G^{pr} \ar@{->}[r]^{\sigma}\ar@{->}[u]^{\overline{\theta}}  & G^{0}
}
\end{equation*}
where $G^{pr}$ is the prefix expansion of $G$. Define the map $\Phi:S \rightarrow G^{0}$ by:
\begin{equation*}
\Phi(s)=\sigma ( m ( \overline{\theta}^{-1} (m(s)))), \Phi(0)=0
\end{equation*}
For each maximal element the preimage under $\overline{\theta}$ is well defined as the map $\theta_{g}$ has the property that $\theta_{g}=\theta_{h} \Rightarrow g=h$ precisely when $\theta_{g} \not = 0 \in S$. Given the preimage is a subset of the F-inverse monoid $G^{pr}$ we know that the maximal element in the preimage is the element $(\lbrace 1,g \rbrace,g)$ for each $g \in G$, from where we can conclude that the map $\sigma$ takes this onto $g \in G$.

We now prove it is a prehomomorphism. Let $\theta_{g},\theta_{h} \in S$, then:
\begin{eqnarray*}
\Phi(\theta_{g})=\sigma ( m(\overline{\theta}^{-1}(\theta_{g}))) = \sigma ( \lbrace 1,g \rbrace, g)= g\\
\Phi(\theta_{h})=\sigma ( m(\overline{\theta}^{-1}(\theta_{h}))) = \sigma ( \lbrace 1,h \rbrace, h)= h\\
\Phi(\theta_{gh})=\sigma ( m(\overline{\theta}^{-1}(\theta_{gh}))) = \sigma ( \lbrace 1,gh \rbrace, gh)= gh
\end{eqnarray*}
Hence whenever $\theta_{g},\theta_{h}$ and $\theta_{gh}$ are defined we know that $\Phi(\theta_{g}\theta_{h})=\Phi(\theta_{g})\Phi(\theta_{h})$. They fail to be defined if:
\begin{enumerate}
\item If $\theta_{gh} = 0$ in $S$ but $\theta_{g}$ and $\theta_{h} \not = 0$ in $S$, then $0=\Phi(\theta_{g}\theta_{h})\leq \Phi(\theta_{g})\Phi(\theta_{h})$

\item If (without loss of generality) $\theta_{g}=0$ then $0=\Phi(0.\theta_{h})= 0.\Phi(\theta_{h})=0$
\end{enumerate}
So prove that the inverse monoid $S$ is strongly 0-F-inverse it is enough to prove then that the map $\Phi$ is idempotent pure, and without loss of generality it is enough to consider maps of only the maximal elements - as the dual prehomomorphism property implies that in studying any word that is non-zero we will be less than some $\theta_{g}$ for some $g \in G$.

So consider the map $\Phi$ applied to a $\theta_{g}$:
\begin{equation*}
\Phi(\theta_{g})=\sigma ( m(\overline{\theta}^{-1}(\theta_{g}))) = \sigma ( \lbrace 1,g \rbrace, g)= g
\end{equation*}

Now assume that $\Phi(\theta_{g}) = e_{G}$. Then it follows that $\sigma (m (\overline{\theta}^{-1}(\theta_{g})))=e_{G}$. As $\sigma$ is idempotent pure, it follows then that $m(\overline{\theta}^{-1}(\theta_{g}))=1$, hence for any preimage $t\in \theta^{-1}(\theta_{g})$ we know that $t \leq 1$, and by the property of being 0-E-unitary it then follows that $t \in E(G^{pr})$. Mapping this back onto $\theta_{g}$ we can conclude that $\theta_{g}$ is idempotent, but by assumption this only occurs if $g = e$.\end{proof}

\begin{proposition}\label{Prop:GrpoidHom}
Let $S = \langle \theta_{g} | g \in G \rangle$ be a strongly 0-F-inverse monoid with maximal elements $Max(S)= \lbrace \theta_{g}:g \in G \rbrace$, where $\theta: G \rightarrow S$ is a dual prehomomorphism. Then the groupoid $G_{\E}$ is Hausdorff, second countable with compact unit space. Also it admits a continuous proper groupoid homomorphism onto the group $G$.
\end{proposition}
\begin{proof}
We record the topological facts about this groupoid here for reference.

Using the map $\Phi$ we construct a map $\rho: \G_{\E} \rightarrow G$ as follows:
\begin{equation*}
\rho([m,\phi]) = \Phi(m)
\end{equation*}
A simple check proves this is a groupoid homomorphism. This map sends units to units as the map $\Phi$ is idempotent pure. We prove continuity by considering preimage of an open set in $G$:
\begin{equation*}
\rho^{-1}(U)=\bigcup_{g \in U}[\theta_{g},\widehat{D}_{\theta^{*}_{g}\theta_{g}}]
\end{equation*}
This is certainly open as each $[\theta_{g},\widehat{D}_{\theta^{*}_{g}\theta_{g}}]$ are elements of the basis of topology of $\G_{\E}$. We check it is proper by observing that for groups $G$ compact sets are finite, and they have preimage:
\begin{equation*}
\rho^{-1}(F)=\bigcup_{g \in F}[\theta_{g},\widehat{D}_{\theta^{*}_{g}\theta_{g}}], \mbox{ } \vert F \vert < \infty 
\end{equation*}
This is certainly compact as these are open and closed sets in the basis of topology for the groupoid $\G_{\E}$.\end{proof}

As $\G_{\widehat{X}} \subseteq \G_{\E}$ we also get a continuous proper groupoid homomorphism from $\G_{\widehat{X}}$ onto a group. We recall a special case of \cite[Lemme 3.12]{MR1703305}.

\begin{lemma}\label{Lem:Lemme}
Let $G$ and $H$ be locally compact, Hausdorff, \'etale topological groupoids and let $\varphi: G \rightarrow H$ be a continuous proper groupoid homomorphism. If $H$ has the Haagerup property then so does $G$. \qed
\end{lemma}

This lets us conclude the following:
\begin{corollary}\label{Cor:Gpoid}
Let $\theta$ be a partial action of $G$ on $X$ such that all the conditions of Proposition \ref{Prop:Aug} are satisfied and such that the inverse monoid $S_{inf}$ is strongly 0-F-inverse. If $G$ has the Haagerup property then so does $\G_{\widehat{X}}$.
\end{corollary}
\begin{proof}
The map induced by the idempotent pure 0-restricted prehomomorphism from $S_{inf}$ to $G$ induces a continuous proper groupoid homomorphism from $G_{\widehat{X}}$ to $G$. This then follows from Lemma \ref{Lem:Lemme}.\end{proof}

\subsection{Partial Actions on Sequences of Graphs}
Let $\lbrace X_{i} \rbrace$ be a sequence of finite graphs with degree $\leq 2k$ and large girth. We begin by considering Pedersens Lemma. The following is \cite[Theorem 6, Chapter XI]{MR1035708}:

\begin{lemma}\label{Lem:GenPet}
Let $X$ be a finite graph. If at most $2k$ edges go into any vertex then all the edges of $X$ can be divided into $k$ classes such that at most two edges from the same class go into any vertex.\qed
\end{lemma}

We remark that we can always assume that the $2k$ here is minimal; there is a smallest even integer that bounds above the degree of all graphs in the sequence. This in particular stops us from doing something unnatural like embedding the 4-regular tree into a 6-regular tree.

\begin{lemma}\label{Lem:SFGL1}
Such a sequence can be almost $k$-oriented and this defines a bounded partial action of $F_{k}$ on $X$
\end{lemma}
\begin{proof}
We work on just the $X_{i}$. Using Lemma \ref{Lem:GenPet}, we partition the edges $E(X_{i})$ into at most $k$ sets $E_{j}$ such that every vertex appears in at most 2 edges from each subset. Pick a generating set $S=\lbrace a_{j} | j\in \lbrace 1,...,k\rbrace\rbrace$ for $F_{k}$ and assign them to the edge sets $E_{j}$, and label the edges that appear in each $E_{j}$ by the corresponding generator. Pedersens Lemma ensures that no more than 2 edges at each point have the same label. This defines a map from the edges to the wedge $\bigvee_{j=i}^{k}S^{1}$. Choose an orientation of each circle and pull this back to the finite graph $X_{i}$ - this provides the partial k-orientation. Now define for each generator the partial bijection $\theta^{i}_{a_{j}}$ that maps any vertex appearing as the source of any edge in $E_{j}$ to the range of that edge. i.e:
\begin{equation*}
\theta^{i}_{a_{j}}(v) := \left\{ \begin{array}{c} r(e) \mbox{ if } \exists e \in E_{j}: s(e)=v \\ \mbox{undefined otherwise} \end{array}\right.
\end{equation*}
For $g=a_{l}^{e_{l}}...a_{m}^{e_{m}}$ we define $\theta^{i}_{g}$ as the product $\theta_{a_{l}^{e_{l}}}...\theta_{a_{m}^{e_{m}}}$; i.e $\theta_{g}^{i}$ moves vertices along any path that is labelled by the word $g$ in the graph $X_{i}$. We observe that for $i \not = i^{'}$ the domain $D_{\theta_{g}^{*}\theta_{g}}^{i}\cap D_{\theta_{g}^{*}\theta_{g}}^{i^{'}}$ is empty hence we can add these partial bijections in $I(X)$ to form $\theta_{g}=\sqcup\theta_{g}^{i}$. It is a remark that as the topology of $X$ is discrete these maps are all continuous and open. It is clear that as each $X_{i}$ is connected that the partial bijections have the property that $\cup_{g} D_{\theta_{g}^{*}\theta_{g}} = X$. Lastly, this map is a dual prehomomorphism as for each $g,h \in G$ we have that $\theta_{g}\theta_{h}=\theta_{gh}$ precisely when both $\theta_{h}$ and $\theta_{gh}$ are defined and moreover if $\theta_{g}\theta_{h}$ is defined then so is $\theta_{gh}$. Hence this collection forms a partial action of $G$ on $X$. We also remark that as each bijection is given translation along a labelling in the free group it is clear that these move elements only a bounded distance and are therefore elements of $I_{b}(X)$.
\end{proof}

We need to show that the partial action generates the metric coarse structure at infinity, we recall the length of a partial bijection:

\begin{definition}
The length of each $\theta_{g}$ is defined to be:
\begin{equation*}
\vert \theta_{g} \vert = \sup \lbrace d(x,\theta_{g}(x)) : x \in Dom(\theta_{g})\rbrace.
\end{equation*} 
\end{definition}

\begin{remark}
As we have a concrete description of each $\theta_{g}$, given on each $X_{i}$, we can see that the length on each $X_{i}$ is given by:
\begin{equation*}
\vert \theta_{g}^{i} \vert = max\lbrace \vert p \vert: p \in \lbrace \mbox{ paths in }X_{i}\mbox{ labelled by g}\rbrace.
\end{equation*} Then $\vert \theta_{g} \vert = \sup_{i} \vert \theta_{g}^{i} \vert$.
\end{remark}

In this situation we require that the partial action contains plenty of infinitely supported elements. 

\begin{proposition}\label{Prop:Inf}
Let $\theta: F_{k} \rightarrow I(X)$ be the dual prehomomorphism corresponding to the bounded partial actions on each $X_{i}$. Then for each $R>0$ there exist finitely many infinite $\theta_{g}$ with $\vert \theta_{g} \vert = \vert g \vert < R$. 
\end{proposition}
\begin{proof}

In the general case we know the following for each $R>0$ and $i \in \mathbb{N}$: $\vert \theta_{g}^{i} \vert \leq \vert g \vert \leq R$. From Lemma \ref{Lem:SFGL1} we know 
that the partial action is defined by moving along paths inside each individual $X_{i}$. So for each $R$ we count the number of words in $F_{k}$ with length less than $R$; this is finite (consider the Cayley graph, which has bounded geometry). Now we observe that on the other hand there are infinitely many simple paths of length less than $R$, thus we must repeat some labellings infinitely many times. These labellings will be contained in words in $F_{k}$ of length less than $R$ hence when we take the supremum we observe that $\vert \theta_{g} \vert = \vert g \vert < R$. \end{proof}

\begin{corollary}
The bounded partial action $\theta$ of $F_{k}$ on $X$ finitely generates the metric coarse structure at infinity, that is the set $\overline{\Delta_{R}}\setminus \Delta_{R} = \bigcup_{\vert g \vert < R}\overline{\Delta_{\theta_{g}}}\setminus \Delta_{\theta_{g}}$ where the index set is finite.
\end{corollary}
\begin{proof}
We proceed by decomposing $\Delta_{R}$ as we did in the proof of Proposition \ref{Prop:Crit}.
\begin{equation*}
\Delta_{R}=(\bigcup_{\vert \theta_{g} \vert < R}\Delta_{\theta_{g}})\cup F_{R}
\end{equation*}
Where $F_{R}$ is the finitely many elements of $\Delta_{R}$ who move between components. We now consider the following decomposition of the set $A:=\lbrace \theta_{g} | \vert \theta_{g}\vert < R\rbrace$ into:
\begin{eqnarray*}
A_{\infty}=\lbrace \theta_{g} | \vert \theta_{g} \vert < R \mbox{ and } \vert supp(\theta_{g})\vert = \infty \rbrace \\
A_{fin}=\lbrace \theta_{g} | \vert \theta_{g} \vert < R\mbox{ and } \vert supp(\theta_{g})\vert < \infty \rbrace
\end{eqnarray*}
The first of these is in bijection with the words in $F_{k}$ that have $\vert g \vert < R$ and define an infinite $\theta_{g}$ from Proposition \ref{Prop:Inf}. Then:
\begin{equation*}
\Delta_{R}=(\bigcup_{g \in A_{\infty}}\Delta_{\theta_{g}})\cup (\bigcup_{g \in A_{fin}}\Delta_{\theta_{g}}) \cup F_{R}
\end{equation*}
We complete the proof by observing that for each $\theta_{g} \in A_{fin}$ the set $\Delta_{\theta_{g}}$ is finite. Therefore:
\begin{equation*}
\overline{\Delta_{R}}\setminus \Delta_{R} = \bigcup_{g \in A_{\infty}}\overline{\Delta_{\theta_{g}}}\setminus \Delta_{\theta_{g}}= \bigcup_{\vert g \vert < R}\overline{\Delta_{\theta_{g}}}\setminus \Delta_{\theta_{g}}
\end{equation*}\end{proof}

\begin{lemma}\label{Lem:ParFree}
The partial action of $F_{k}$ defined above extends to $\beta X$ and is free on the boundary $\partial \beta X$.
\end{lemma}
\begin{proof}
Fix $g \in F_{k}$. As the sequence of graphs has large girth there exists an $i_{g} \in \mathbb{N}$ such that for all $i\geq i_{g}$ $\theta_{g}^{i}$ has no fixed points.  Let $\omega \in \widehat{D}_{\theta_{g}^{*}\theta_{g}}$ and assume for a contradiction that $\theta_{g}(\omega)=\omega$. Consider the graph $G$ with vertex set $V=\sqcup_{i\geq i_{g}} X_{i}$, where two vertices $x$ and $y$ are joined by an edge if and only if $y=\theta_{g}(x)$. This graph will have degree at most 2. Additionally, the set of non-isolated vertices inside $V$ contains $D_{\theta_{g}^{*}\theta_{g}}$ and so the subgraph with vertex set consisting of the non-isolated vertices and the same edge set is chosen by $\omega$ and has degree at most 2. Any such graph can be at most $3$ coloured, subsequently the vertex set breaks into three disjoint pieces $V_{1}, V_{2}$ and $V_{3}$ and the ultrafilter $\omega$ will pick precisely one of the $V_{i}$. The action of $\theta_{g}$ sends $V_{i}$ into $V_{i+1}$ modulo $3$ by construction. Lastly, $V_{i} \in \omega$ implies $V_{i+1} \in \theta_{g}(\omega)=\omega$, which is a contradiction.
\end{proof}

This freeness gives us a tool to understand the structure of $S_{\inf}$.

\begin{lemma}\label{Lem:CP}
Let $\lbrace X_{i} \rbrace$ be a sequence of graphs and let $G$ be a group which acts partially on each $X_{i}$. If $G$ fixes any sequence in $\lbrace X_{i} \rbrace$ then the partial action is not free on $\partial \beta X$. 
\end{lemma}
\begin{proof}Let $\theta_{g}$ denote the disjoint union of the $\theta_{g}^{i}$ arising from the partial action of $G$ on each $X_{i}$. To prove this it is enough to show that there is a single $\omega \in\partial\beta X$ that is fixed by the action of some $g \in G$. The hypothesis that $G$ fixes a sequence gives us $\textbf{x}:=\lbrace x_{n} \rbrace_{I}$ with $I$ infinite and $\theta_{g}(\lbrace x_{n} \rbrace) = \lbrace\theta_{g}^{n}(x_{n}) \rbrace_{I} =\lbrace x_{n} \rbrace_{I}$.

Now consider an ultrafilter $\omega \in \partial\beta X$ that picks $\textbf{x}$. Then this ultrafilter $\omega$ is an element of $D_{g^{*}g}$ as $\textbf{x} \subset D_{g^{*}g}$.  Now for any $A \in \omega$ and consider the intersection $A \cap \textbf{x}$. This is fixed by the action of $g$, as it is a subset of $\textbf{x}$. Hence we have: $\theta_{g}(A \cap \textbf{x}) \in \theta_{g}(\omega)$ for every $A \in \omega$. As $\theta_{g}(\omega)$ is an ultrafilter $A \in \theta_{g}(\omega)$, so in particular $\omega \subseteq \theta_{g}(\omega)$, whence $\theta_{g}(\omega)=\omega$. 
\end{proof}

Recall that the inverse monoid $S_{inf}$ is represented geometrically by partial bijections on $I(X)$. This representation gives us access to the geometry of $X$, which we can utilise, in addition to Lemma \ref{Lem:CP}, to understand the structure of $S_{inf}$. 

\begin{lemma}
Consider the inverse monoid $S_{inf}$ as a submonoid of $I(X)$. Then the following hold:
\begin{enumerate}
\item $S_{inf}$ has the property that $(g\not = e_{G}$ and $\theta_{g}\not = 0)$ implies $\theta_{g}$ is not an idempotent; 
\item $S_{inf}$ is $0$-E-unitary;
\item $S_{inf}$ has maximal element set $\lbrace \theta_{g} : g \in F_{k} \rbrace$.
\end{enumerate}
\end{lemma}
\begin{proof}
\hskip 1pt

\begin{enumerate}
\item We prove that no non-zero $\theta_{g}$ are idempotent. To do this we pass to the induced action on $\beta X$. We observe that if $\theta_{g}$ is idempotent on $X$ then it extends to an idempotent on $\beta X$, hence on the boundary $\partial \beta X$. $\theta_{g}$ is non-zero implies that there is a non-principal ultrafilter $\omega$ in the domain $\widehat{D}_{\theta_{g}}$. The result then follows from the observation that $\theta_{g}\circ \theta_{g} (\omega) = \theta_{g}(\omega)$ implies that $\theta_{g}$ must now fix the ultrafilter $\theta_{g}(\omega)$, which by Lemma \ref{Lem:ParFree} cannot happen.

\item For 0-E-unitary it is enough to prove that $f \leq \theta_{g}$ implies $\theta_{g} \in E(S)$. Again, we extend the action to $\beta X$. We observe that if $\theta_{g}$ contains an idempotent, then we can build a sequence of elements of $x_{i} \in f \cap \widehat{D}_{\theta_{g}^{*}\theta_{g}} \cap X_{i}$ such that $\theta_{g}$ fixes the sequence, and hence fixes any ultrafilter $\omega$ that picks this sequence by Lemma \ref{Lem:CP}. This is a contradiction, from where we deduce that the only situation for which $f \leq \theta_{g}$ is precisely when $g =e_{G}$ hence trivially $e \leq \theta_{g}$ implies $\theta_{g} \in E(S)$. For the general case, we remark that by the above statement coupled with the dual prehomomorphism property shows that $f \leq s$ implies $s \leq \theta_{e_{G}}$, hence is an idempotent.  

\item We construct the maximal elements. Observe that using the dual prehomomorphism it is clear that every non-zero word $s \in S$ lives below a non-zero $\theta_{g}$. So it is enough to prove that for $\theta_{g}, \theta_{h} \not = 0$, $\theta_{g} \leq \theta_{h} \Rightarrow \theta_{g}=\theta_{h}$. Let $\theta_{g} \leq \theta_{h}$. This translates to $\theta_{h}\theta_{g}^{*}\theta_{g} = \theta_{g}$, hence for all $x \in \widehat{D}_{\theta_{g}^{*}\theta_{g}}: \theta_{h}(x) = \theta_{g}(x)$. Hence $\theta_{g}^{*}\theta_{h} \in E(S)$. From here we see that $\theta_{g}^{*}\theta_{h} \leq \theta_{e}$. From (2) we can deduce: $\theta_{g}^{*}\theta_{h}\leq \theta_{g^{-1}h}$ implies $\theta_{g^{-1}h} \in E(S_{inf})$. By (1) this implies $\theta_{g^{-1}h}=\theta_{e}$, and this happens if and only if $g^{-1}h = e$, i.e $g=h$.
\end{enumerate}
\end{proof}
Appealing to the machinery we developed earlier in Propositions \ref{Prop:Strongly} and \ref{Prop:GrpoidHom} we get the following corollary immediately.

\begin{corollary}\label{Cor:MT}
The inverse monoid $S_{inf}$ is strongly 0-F-inverse.
\end{corollary}

We now have enough tools to prove the following result:
\begin{theorem}
Let $\lbrace X_{i} \rbrace$ be a sequence of finite graphs of large girth and vertex degree uniformly bounded above by $2k$ and let $X$ be the corresponding space of graphs. Then the boundary coarse Baum-Connes conjecture holds for $X$. 
\end{theorem}
\begin{proof}
Combining the results in the previous section we know that a free group of rank $k$ acts partially freely on the boundary $\partial\beta X$ in such a way as to give us a second countable locally compact Hausdorff \'etale topological groupoid $\G_{\widehat{X}}$. This groupoid implements an isomorphism $G(X)|_{\partial\beta X} \cong \partial\beta X \rtimes \G_{\widehat{X}}$. We also know using Corollary \ref{Cor:MT} that the inverse monoid generated by the infinite support elements $S_{inf}$ is a strongly $0$-F-inverse monoid, admitting a $0$-restricted idempotent pure prehomomorphism onto $F_{k}$. Hence the groupoid $\G_{\widehat{X}}$ admits a proper continuous groupoid homomorphism onto $F_{k}$, and so has the Haagerup property by Corollary \ref{Cor:Gpoid}. 

We can now conclude the Theorem by remarking that the isomorphism of groupoids $G(X)|_{\partial\beta X} \\ \cong \partial\beta X \rtimes \G_{\widehat{X}}$ turns the Baum-Connes conjecture for $G(X)|_{\partial\beta X}$ into a specific case of the Baum-Connes conjecture for $\G_{\widehat{X}}$. As $\G_{\widehat{X}}$ has the Haagerup property we can conclude that the Baum-Connes assembly map with any coefficients is an isomorphism \cite{MR1798599} and so the assembly map required for the boundary conjecture is also an isomorphism.
\end{proof}

\section{Applications to the coarse Baum-Connes conjecture.}\label{Sect:apps}
In this section we capture, via homological methods, all the known counterexample arguments present in the literature using the boundary conjecture. We begin by confirming that the kernel $I$ constructed in section \ref{Sect:CE} is the ghost ideal $I_{G}$. In order to prove this we need some technology of \cite{MR1911663}:

\begin{lemma}\label{Lem:Lemma9}\cite[Lemma 9]{MR1911663}
If an \'etale topological groupoid $\G$ acts on a $C^{*}$-algebra A, then the map $C_{c}(\G,A) \rightarrow C_{0}(\G,A)$ extends to an injection (functorial in A) from $A\rtimes_{r} \G$ to $C_{0}(\G,A)$. \qed
\end{lemma}

\begin{remark}
The phrase ``functorial in $A$'' allows us, given a map: $A \rightarrow B$ of $\G-C^{*}$-algebras, to build the following square:
\begin{equation*}
\xymatrix{A\rtimes_{r} \G \ar@{->}[r] \ar@{^{(}->}[d] &B \rtimes_{r} \G \ar@{^{(}->}[d]\\
C_{0}(\G ,A)\ar@{->}[r] &   C_{0}(\G ,B)
}
\end{equation*}
\end{remark}

\begin{remark}
The map provided is not a $*-$homomorphism as it takes convolution in $A\rtimes_{r} \G$ to point wise multiplication in $C_{0}(\G,A)$. However it suffices for applications as the map is continuous.
\end{remark}

\begin{proposition}
Let $X$ be a uniformly discrete bounded geometry metric space. Then the kernel of the map
\begin{equation*}
l^{\infty}(X,\mathcal{K})\rtimes_{r}G(X) \rightarrow (l^{\infty}(X,\mathcal{K})/C_{0}(X,\mathcal{K}))\rtimes_{r}G(X)|_{\partial\beta X}
\end{equation*}
is the ghost ideal $I_{G}$.
\end{proposition}
\begin{proof}
Lemma \ref{Lem:Lemma9} implies that the following diagram commutes:
\begin{equation*}
\xymatrix{\ell^{\infty}(X,\mathcal{K})\rtimes G(X) \ar@{->}[r]^{q} \ar@{^{(}->}[d]^{i^{'}} &\frac{\ell^{\infty}(X,\mathcal{K})}{C_{0}(X,\mathcal{K})}\rtimes G(X) \ar@{^{(}->}[d]^{i}\\
C_{0}(G(X),\ell^{\infty}(X,\mathcal{K}))\ar@{->}[r]^{q^{'}} &   C_{0}(G(X),\frac{\ell^{\infty}(X,\mathcal{K})}{C_{0}(X,\mathcal{K})})
}
\end{equation*}
The downward maps being injective implies that the kernel is precisely the kernel of induced map:
\begin{equation*}
q:C^{*}X=\ell^{\infty}(X,\mathcal{K})\rtimes_{r}G(X) \rightarrow C_{0}(G(X),l^{\infty}(X,\mathcal{K})/C_{0}(X,\mathcal{K})).
\end{equation*}
Using the definitions we can compute this kernel:
\begin{eqnarray*}
I & =&  \lbrace f | i^{'}(f) \in C_{0}(G(X),C_{0}(X,\mathcal{K}))\rbrace \\
& = & \lbrace f | \forall \epsilon > 0 \exists K \subset X \times X \mbox{ compact }: \vert f_{xy} \vert \leq \epsilon \forall (x,y) \in X\times X \setminus K \rbrace.
\end{eqnarray*}
As $X$ is uniformly discrete with bounded geometry and $X \times X$ is equipped with the product topology we can replace compact by bounded, which using Definition \ref{Def:Ghost} is the ghost ideal $I_{G}$.
\end{proof}
Recall that the assembly map $\mu_{I_{G}}$ associated to the open saturated subset $X$ is given by the composition of $\mu_{X \times X}: K^{top}_{*}(X\times X) \rightarrow K_{*}(\mathcal{K}(\ell^{2}(X))$ with the inclusion $i_{*}:K_{*}(\mathcal{K}(\ell^{2}(X)) \rightarrow K_{*}(I_{G})$. So to understand how the assembly map $\mu_{I_{G}}$ behaves, it is enough to consider the behaviour of the inclusion $i_{*}$ as the map $\mu_{X \times X}$ is an isomorphism.

\begin{proposition}\label{Prop:Ghost}
Let $\lbrace X_{i} \rbrace_{i\in \mathbb{N}}$ be a sequence of finite graphs. The the following hold for the space of graphs $X$:
\begin{enumerate}
\item The induced map $i_{*}:K_{*}(\mathcal{K}(\ell^{2}(X,\mathcal{K})) \rightarrow K_{*}(C^{*}X)$ is injective.
\item If $X$ is an expander then $K_{*}^{top}(X\times X)\cong K_{*}(\mathcal{K}(\ell^{2}(X)) \rightarrow K_{*}(I_{G})$ is not surjective.
\end{enumerate}
\end{proposition}
\begin{proof}
Consider the coarse map $X=\sqcup_{i} X_{i} \rightarrow P:=\sqcup_{i} \ast$ given by projecting each factor to a point and $P$ carries the coarse disjoint union metric. This induces a tracelike map on the K-theory of $C^{*}X$, which we denote by $Tr$. Observe also that $C^{*}P \cong \ell^{\infty}(P,\mathcal{K})+\mathcal{K}(\ell^{2}(P,\mathcal{K}))$, hence has $K_{0}(C^{*}P) \cong  \frac{\prod\mathbb{Z}}{K}$, where $K$ is the subgroup of $\bigoplus \mathbb{Z}$ of elements with sum equal to 0. Now:

To prove (1), Observe that the following diagram commutes:
\begin{equation*}
\xymatrix{ K_{0}(C^{*}X) \ar[r]& \frac{\prod\mathbb{Z}}{K}\\
K_{0}(\mathcal{K}) \ar[u] \ar[ur] &
}
\end{equation*}
Where the image of a rank one projection in $K_{0}(\mathcal{K})$ is given by a vector $(1,0,...)$, which is certainly an injection. It now follows that the upward map is injective.

For (2) consider a nontrivial ghost projection $p \in I_{G}$. It is well known that $Tr([p]) \not\in \frac{\bigoplus \mathbb{Z}}{K}$, whilst $Tr([k]) \in \frac{\bigoplus\mathbb{Z}}{K}$ for any compact operator \cite{explg1}. As these differ under $Tr$, they cannot possibly be equal in $K_{0}(I_{G})$.
\end{proof}

The boundary coarse Baum-Connes conjecture has applications to the coarse Baum-Connes conjecture for spaces of graphs:

\begin{theorem}\label{Prop:Cor}
Let $X$ be a space of graphs arising from a sequence of finite graphs $\lbrace X_{i} \rbrace$. If  boundary coarse Baum-Connes conjecture is injective then the following hold:
\begin{enumerate}
\item The coarse Baum-Connes assembly map for $X$ is injective.
\item If $X$ is an expander then the coarse Baum-Connes assembly map for $X$ fails to be surjective.
\end{enumerate}
\end{theorem}
\begin{proof}
Consider the diagram:
$$
\xymatrix@=1em{& & & \frac{\prod_{i\in\mathbb{N}}K_{0}(C^{*}X_{i})}{\bigoplus_{i\in\mathbb{N}}K_{0}(C^{*}X_{i})}& & &\\
\ar[r] & K_{1}(C^{*}_{r}(G(X)|_{\partial\beta X}) \ar[r] & K_{0}(I_{G})\ar[ur] \ar[r]& K_{0}(C^{*}_{r}(G(X)))\ar[u]^{d_{*}} \ar[r]& K_{0}(C^{*}_{r}(G(X)|_{\partial\beta X})\ar[r] & K_{1}(I_{G}) \ar[r] & \\
\ar[r] & K_{1}^{top}(G(X)|_{\partial\beta X}) \ar[r] \ar@{^{(}->}[u]& \mathbb{Z} \ar@{^{(}->}[ur] \ar[r]\ar@{^{(}->}[u]& K_{0}^{top}(G(X)) \ar[r]\ar[u]^{\mu}& K_{0}^{top}(G(X)|_{\partial\beta X}) \ar[r]\ar@{^{(}->}[u]& 0\ar[u]^{0} \ar[r] & 
}
$$
We prove (1) by considering an element $x \in K_{0}^{top}(G(X))$ such that $\mu(x)=0$. Then $x$ maps to $0\in K_{0}^{top}(G(X)|_{\partial\beta X})$ and so comes from an element $y \in \mathbb{Z}$. Each square commutes hence $y$ maps to $0 \in K_{0}(C^{*}_{r}(G(X))$. As the composition up and left (as indicated in the diagram) is injective by Proposition \ref{Prop:Ghost}, we know that $y \in \mathbb{Z}$  is in fact $0 \in \mathbb{Z}$. Hence $x=0$.

To see (2): take any non-compact ghost projection $p \in K_{0}(I_{G})$, which does not lie in the image of $\mathbb{Z}$ as it does not vanish under the trace $d_{*}$. Push this element to $q \in K_{0}(C^{*}_{r}(G(X))$. Assume for a contradiction that $\mu$ is surjective. Then there is an element $x \in \mathbb{Z}$ that maps to $q$, and so $d_{*}(q)=0$, as the image of compact operators lies in the kernel of the map $d_{*}$. However, we know that $d_{*}(q)=d_{*}(p)$ is certainly non-zero.
\end{proof}

\subsection{Some Remarks on the Max Conjecture.}
\begin{proposition}\label{Prop:Max}
Let $X$ be the space of graphs arising from a sequence of finite graphs $\lbrace X_{i} \rbrace$. Then
\begin{enumerate}
\item the maximal coarse Baum-Connes assembly map is an isomorphism if and only if the maximal Boundary coarse Baum-Connes map is an isomorphism.
\item the maximal coarse assembly map is injective if and only if the maximal boundary assembly map is injective.
\end{enumerate}
\end{proposition}
\begin{proof}
As before we consider a diagram, this time of maximal algebras:
\begin{equation*}
\xymatrix@=0.7em{
K_{1}(C^{*}(G(X)|_{\partial\beta X}) \ar[r] & K_{0}(\mathcal{K}) \ar[r]& K_{0}(C^{*}(G(X))) \ar[r]& K_{0}(C^{*}(G(X)|_{\partial\beta X})\ar[r] & K_{1}(\mathcal{K})  \\
K_{1}^{top}(G(X)|_{\partial\beta X}) \ar[r] \ar[u]& K_{0}^{top}(X \times X) \ar[r]\ar[u]^{\ucong}& K_{0}^{top}(G(X)) \ar[r]\ar[u]& K_{0}^{top}(G(X)|_{\partial\beta X}) \ar[r]\ar[u]^{\mu_{bdry}}& K_{1}^{top}(X \times X)\ar[u]^{\ucong}
}
\end{equation*}
Both parts follow from a diagram chase.
\end{proof}

This result captures completely \cite[Corollary 4.18]{MR2568691}, but the proof is very much more elementary. It is also clear that the argument above also works in the case of injectivity, which is related to the content of \cite[Theorem 5.1]{MR2431253}.

\section{A Counterexample to the Boundary Conjecture}\label{Sect:Counter}
In this section we develop the ideas of Higson, Lafforgue and Skandalis concerning the counterexamples to the coarse Baum-Connes conjecture further, to construct a space of graphs $Y$ that has exceptional properties at infinity. The main idea is to decompose the boundary groupoid further, giving a new short exact sequence at infinity similar to the sequences considered in Chapter 4. From this, we then construct an operator that is not a ghost operator, but is ghostly on certain parts of the boundary. A tracelike argument, similar to those of \cite{higsonpreprint, explg1} then allows us to conclude that the boundary coarse Baum-Connes conjecture fails to be surjective for the space $Y$.

\subsection{The space and its non-ghosts.}

The space we are going to consider first appeared in \cite{MR2363697}.

Let $\lbrace X_{i} \rbrace_{i \in \mathbb{N}}$ be a sequence of finite graphs. Then we construct a space of graphs in the following manner: Let $Y_{i,j} = X_{i}$ for all $j \in \mathbb{N}$ and consider $Y:= \sqcup_{i,j} Y_{i,j}$. We metrize this space using a box metric - that is with the property that $d(Y_{i,j},Y_{k,l}) \rightarrow \infty$ as $i+j+k+l \rightarrow \infty$. 

Now let $\lbrace X_{i} \rbrace_{i}$ be an expander sequence. As discussed in Section \ref{Sect:GO} of Chapter 4, we can construct a ghost operator $p= \prod_{i} p_{i}$ on $X$, the space of graphs of $\lbrace X_{i} \rbrace_{i}$. Similarly, we can construct this operator on $Y$. In this situation we get a projection $q:=\prod_{i,j}p_{i} \in C^{*}_{u}Y$, which is a constant operator in the $j$ direction. This was precisely the operator of interest in \cite{MR2363697}, as it can be seen that $q$ is not a ghost operator, as its matrix entries do not vanish in the $j$ direction - a fact proved below in Lemma \ref{Lem:nag}.

Recall that associated to $Y$ we have a short exact sequence of $C^{*}$-algebras:
\begin{equation*}
\xymatrix{
0 \ar[r] & ker(\pi) \ar[r]& C_{r}^{*}(G(Y)) \ar[r]^{\pi} & C_{r}^{*}(G(Y)|_{\partial\beta Y}) \ar[r] & 0.
}
\end{equation*}

We remark the kernel, $ker(\pi)$ consists of all the ghost operators in $C^{*}_{u}(Y)$, that is those operators with matrix coefficients that tend to $0$ in all directions on the boundary. 

\begin{lemma}\label{Lem:nag}
The projection $\pi(q):= \pi(\prod_{i,j}p_{i}) \not = 0 \in C^{*}_{r}(G(Y)|_{\partial\beta Y})$. That is $q \not\in ker(\pi)$.
\end{lemma}
\begin{proof}
We first observe that every bounded subset $B$ of $Y$ is contained in some rectangle of the form $R_{i_{B},j_{B}}:=\sqcup_{i\leq i_{B},j\leq j_{B}}Y_{i,j}$. So to prove that $q$ is not a ghost operator it suffices to show that there exists an $\epsilon>0$ such that for all rectangles $R_{i,j}$ there is a pair of points $x,y$ in the compliment of the rectangle such that  the norm $\Vert q_{x,y} \Vert \geq \epsilon$. To prove this, recall that the projection $q$ is a product of projections $p_{i}$ on each $X_{i}$ and fixing $j$, these projections form a ghost operator. 

Fix $\epsilon = \frac{1}{2}$. Then there exists an $i_{\epsilon}$ with the property that $\forall i>i_{\epsilon}$ and for every $x,y \in \sqcup_{i}X_{i}$ we know $\Vert p_{i,x,y} \Vert < \epsilon$. We remark that this $i_{\epsilon}$ can be taken to be the smallest such. So for $i \leq i_{\epsilon}-1$, we have that $\Vert p_{i,x,y} \Vert \geq \frac{1}{2}$. Now let $R_{i_{\epsilon}-1,\infty}$ be the vertical rectangle $\sqcup_{i\leq i_{\epsilon}-1,j} Y_{i,j}$. 

To finish the proof, consider an arbitrary finite rectangle $R_{i,j}$. This intersects the infinite rectangle $R_{i_{\epsilon}-1,\infty}$ in a bounded piece. Now pick any pair of points in a fixed box $x,y \in Y_{k,l} \subset R_{i_{\epsilon}-1,\infty} \setminus R_{i,j}$. Then for those points $x,y$ it is clear that $\Vert q_{x,y} \Vert = \Vert p_{k,x,y}\Vert \geq \frac{1}{2}$.
\end{proof}

We now describe the boundary $\partial\beta Y$. We are aiming at a decomposition into saturated pieces and with that in mind we construct a map to $\beta X$.

Consider the map $\beta Y \twoheadrightarrow \beta X \times \beta \mathbb{N}$ induced by the bijection of $Y$ with $X \times \mathbb{N}$ and the universal property of $\beta Y$. Now define:
\begin{equation*}
f: \beta Y \rightarrow \beta X \times \beta \mathbb{N} \rightarrow \beta X
\end{equation*}
The map $f$ is continuous, hence the preimage of $X$ under projection onto the first factor is an open subset of $\beta Y$, which intersects the boundary $\partial \beta Y$. In fact, what we can see is that $f^{-1}(X)= \sqcup f^{-1}(X_{i})$, where each $f^{-1}(X_{i})$ is closed, and therefore homeomorphic to $X_{i} \times \beta \mathbb{N}$. We can define $U = f^{-1}(X)\cap \partial\beta Y$.

\subsection{The boundary groupoid associated to the box space of a discrete group with the Haagerup property.}

Let $\Gamma$ be a finitely generated residually finite discrete group with the Haagerup property, and let $\lbrace N_{i}\rbrace$ be a family of nested finite index subgroups with trivial intersection. Let $X_{i}:=Cay(\frac{\Gamma}{N_{i}})$. In this context, the boundary groupoid is generated by the action of the group $\Gamma$ extended to the boundary (see Proposition \ref{Prop:Crit}). In this context we can show that $U$ defined above is saturated:

\begin{lemma}
$U$ is an open, saturated subset of the boundary $\partial\beta Y$. 
\end{lemma}
\begin{proof}
We have already shown above that $U$ is open. To see it is saturated we prove that $U^{c}$ is saturated. Let $g_{Y}$ and $g_{X}$ be the continuous extensions of the map obtained by acting through $g$ on $Y$ and $X$ respectively. We observe that the following diagram commutes:
\begin{equation*}
\xymatrix{
\overline{g}_{Y}:\beta Y\ar[r]\ar[d]^{p} & \beta Y\ar[d]^{p}\\
\overline{g_{X} \times 1}:  \beta X \times \beta \mathbb{N} \ar[r] & \beta X \times \beta \mathbb{N}
}
\end{equation*}
The projection onto $\beta X$ is also equivariant under this action. Assume for a contradiction that $U^{c}$ is not saturated; there exists $\gamma$ in $U^{c}$ such that $\overline{g}_{Y}(\gamma) \in U$. It follows that $\overline{g_{X} \times 1}(p(\gamma))$ is in $p(U)$, whilst $p(\gamma) \in p(U^{c})$, hence $\overline{g_{X}}(f(\gamma))\in U$ whilst $f(\gamma) \in U^{c}$. This is a contradiction as $f(U^{c}) = \partial\beta X$ is saturated.
\end{proof}

This gives us two natural complimentary restrictions of $G(Y)|_{\partial\beta Y}$ and a short exact sequence of function algebras as in Chapter 4:
\begin{equation*}
\xymatrix{
0 \ar[r] & C_{c}(G(Y)|_{U}) \ar[r]& C_{c}(G(Y)|_{\partial\beta Y}) \ar[r] & C_{c}(G(Y)|_{U^{c}}) \ar[r] & 0.
}
\end{equation*}

We will now show that the corresponding sequence:
\begin{equation*}
\xymatrix{
0 \ar[r] & C^{*}_{r}(G(Y)|_{U}) \ar[r]& C_{r}^{*}(G(Y)|_{\partial\beta Y}) \ar[r]^{h} & C_{r}^{*}(G(Y)|_{U^{c}}) \ar[r] & 0
}
\end{equation*}
fails to be exact in the middle. We proceed as in \cite{explg1,MR1911663} by using the element $\pi(q)$, which certainly vanishes under the quotient map from $C^{*}_{r}(G(Y)|_{\partial\beta Y}) \rightarrow C^{*}_{r}(G(Y)|_{U^{c}})$. To show the failure we will show this sequence fails to be exact in the middle at the level of K-theory and for this we will require a firm understanding of the structure of $G(Y)|_{U}$.

We observe the following facts: 
\begin{enumerate}
\item $\Gamma$ acts on the space $Y:=\sqcup_{i,j}Y_{i,j}$ built from $\lbrace X_{i} \rbrace$.
\item This action becomes free on piece of the boundary that arises as $i \rightarrow \infty$, that is $\Gamma$ acts freely on $U^{c}$.
\item The group action generates the metric coarse structure on the boundary; the finite sets associated to each $R>0$ in the decomposition are now finite rectangles. This follows from considerations of the metric on $Y$.
\end{enumerate} 

It follows from the proof of Proposition \ref{Prop:Crit} that the groupoid $G(Y)|_{U^{c}}$ is isomorphic to $U^{c}\rtimes \Gamma$ and under the assumption that $\Gamma$ has the Haagerup property we can conclude that the Baum-Connes assembly map for the groupoid $G(Y)|_{U^{c}}$ is an isomorphism (with any coefficients). We now concern ourselves with $G(Y)|_{U}$.

\begin{lemma}\label{Lem:CE3}
The groupoid $G(Y)|_{U}$ is isomorphic to $\sqcup_{i}(X_{i}\times X_{i})\times G(\mathbb{N})|_{\partial\beta \mathbb{N}}$.
\end{lemma}
\begin{proof}
We consider the preimages $f^{-1}(X_{i})$. These are clearly invariant subsets of $\beta Y$ that when intersected with the boundary $\partial\beta Y$ are contained within $U$. The restriction of \\$(G(Y)|_{\partial\beta Y})|_{f^{-1}(X_{i})}$ for each $i$ isomorphic to the closed subgroupoid $G(X_{i}\times \mathbb{N})|_{\partial\beta \mathbb{N}}$ of $G(Y)|_{U}$. These groupoids are disjoint by construction and therefore the inclusion $\sqcup_{i}G(X_{i}\times \mathbb{N})|_{\partial\beta \mathbb{N}}$ is an open subgroupoid of $G(Y)|_{U}$. We now prove that:
\begin{enumerate}
\item each $G(X_{i} \times \mathbb{N})|_{\partial\beta \mathbb{N}}$ is isomorphic to $(X_{i}\times X_{i})\times G(\mathbb{N})|_{\partial\beta \mathbb{N}}$, where $\mathbb{N}$ has the well-spaced metric;
\item the union $\sqcup_{i}G(X_{i}\times \mathbb{N})|_{\partial\beta \mathbb{N}}$ is the entire of $G(Y)|_{U}$.
\end{enumerate}
To prove (1), observe that the groupoid decomposes as $G(X_{i}\times \mathbb{N}) = \bigcup_{R>0}\overline{\Delta_{R}(X_{i}\times \mathbb{N})}$. For each $R>0$ we can find a $j_{R}$ such that $\Delta_{R}(X_{i}\times \mathbb{N}) = F_{R} \cup \bigcup_{j>j_{R}}\Delta_{R}(X_{i}\times \lbrace j \rbrace)$, hence for the boundary part of this groupoid it is enough to understand what happens in each piece $Y_{i,j}$, which is constant for each $j$. Secondly, observe that in the induced metric on a column, the pieces $Y_{i,j}$ separate as $j\rightarrow \infty$. This, coupled with the fact that for large enough $R$, we know that $ \Delta^{j}_{R}(X_{i}\times \mathbb{N}) = X_{i} \times X_{i}$ allow us to deduce that any behaviour at infinity of this groupoid is a product of $X_{i} \times X_{i}$ and the boundary groupoid $G(\mathbb{N})|_{\partial\beta \mathbb{N}}$ where $\mathbb{N}$ has the coarsely disconnected metric. This groupoid is isomorphic to $\partial\beta \mathbb{N}$, from which we can deduce that $G(X_{i}\times \mathbb{N})|_{\partial\beta \mathbb{N}} = (X_{i}\times X_{i})\times \partial\beta \mathbb{N}$ for each $i$.

To prove (2) we assume for a contradiction that there is a partial translation $t$, such that $\overline{t}$ is not an element of the disjoint union. Such an element maps some $(x_{i}, \omega)$ to $(x_{k},\omega)$, where $i\not =k$. Without loss of generality assume also $t$ has translation length at most $R$. Then the domain and range of $t$ are both infinite (as the closure is defined in $G(Y)|_{U}$), and must be contained within a strip of width at most $R>0$. From the definition of the metric, there are only finitely many $Y_{i,j}$ within such a rectangle, hence $t \in F_{R}$ and hence $\overline{t}$ is not defined in $G(Y)|_{U}$, which yields a contradiction.
\end{proof}

\begin{remark}
Lemma \ref{Lem:CE3} allows us to conclude that $C^{*}_{r}(G(Y)|_{U}) \cong \bigoplus_{i}M_{\vert X_{i}\vert} \otimes C(\partial\beta Y)$
\end{remark}

To conclude that $[\pi(q)]$ is not an element of $K_{0}(C^{*}_{r}(G(Y)|_{U}))$ we construct a trace-like map.

\begin{theorem}
The element $\pi(q)$ maps to $0$ in $C^{*}_{r}(G(Y)|_{U^{c}})$, but does not belong to the image of  $K_{0}( C^{*}_{r}(G(Y)|_{U}))$ in $K_{0}(C^{*}_{r}(G(Y)|_{\partial\beta Y}$.
\end{theorem}
\begin{proof}
The first part follows from the definition of $\pi(q)$; the quotient map $h$ kills all operators that are ghostly in the direction $i \rightarrow \infty$ and $\pi(q)$ is such an operator.

To prove the second component we remark that each $U_{i}:=f^{-1}(X_{i}) \cap U$ is a closed saturated subset of $\partial\beta Y$, hence we can consider the reduction to $U_{i}$ for each $i$. We consider the product, and the following map:
\begin{eqnarray*}
\phi : C^{*}_{r}(G(Y)|_{\partial\beta Y})& \rightarrow & \prod_{i} C^{*}_{r}(G(Y)|_{U_{i}})= \prod_{i}C^{*}_{r}(G(X_{i} \times \mathbb{N})\\
 T  &\mapsto & \prod_{i}T|_{U_{i}}
\end{eqnarray*}
Under the map $\phi$, the ideal $C^{*}_{r}(G(Y)|_{U})= \bigoplus_{i}M_{\vert X_{i}\vert} \otimes C(\partial\beta Y)$ maps to the ideal \\$\bigoplus_{i}C^{*}_{r}(G(X_{i} \times \mathbb{N}))$. So, we can define a tracelike map, in analogy to \cite[Section 6]{explg1}, by composing with the quotient map $\tau$ onto $\frac{\prod_{i}C^{*}_{r}(G(X_{i} \times \mathbb{N}))}{\bigoplus_{i}C^{*}_{r}(G(X_{i} \times \mathbb{N}))}$. This gives us a map at the level of K-theory:
\begin{equation*}
Tr_{*}=\phi \circ \tau : K_{0}(C^{*}_{r}(G(Y)|_{\partial\beta Y})) \rightarrow \frac{\prod_{i}K_{0}(C^{*}_{r}(G(X_{i} \times \mathbb{N})))}{\bigoplus_{i}K_{0}(C^{*}_{r}(G(X_{i} \times \mathbb{N})))}= \frac{\prod_{i}K_{0}(C(\partial\beta \mathbb{N}))}{\bigoplus_{i}K_{0}(C(\partial\beta \mathbb{N}))}
\end{equation*}

By construction, $K_{0}(C^{*}_{r}(G(Y)|_{U}))$ vanishes under $Tr_{*}$. We now consider $[\pi(q)]$ under $Tr_{*}$. Recall that $q=\prod_{i,j}p_{i}$. We define $q_{i}=\prod_{j}p_{i}$ and observe that the operation of reducing to $G(Y)|_{U_{i}}$ can be performed in two commuting ways: restricting to $U$ then $f^{-1}(X_{i})$ or by restricting to $f^{-1}(X_{i})$ then $U$. The second tells us that $q_{i}=p_{i} \otimes 1_{\beta\mathbb{N}}$ is constant in the $j$ direction and when restricted to the boundary is $\pi(q_{i})=p_{i}\otimes 1_{\partial\beta \mathbb{N}}$. Hence, $Tr_{*}([\pi(q)])=[1_{\partial\beta \mathbb{N}},1_{\partial\beta\mathbb{N}},...] \not = 0$ and so $[\pi(q)] \not \in K_{0}(C^{*}_{r}(G(Y)|_{U})$.
\end{proof}

So in this case we have the following diagram:
\begin{equation*}
\xymatrix@=0.7em{
 K_{1}(C(U^{c})\rtimes \Gamma) \ar[r] & K_{0}(\ker (\pi)) \ar[r]& K_{0}(C^{*}_{r}(G(Y)|_{\partial\beta Y})) \ar[r]& K_{0}(C(U^{c})\rtimes \Gamma)\ar[r] & K_{1}(\ker (\pi))  \\
 K_{1}^{top}(U^{c}\rtimes \Gamma) \ar[r] \ar[u]^{\ucong}& K_{0}^{top}(G(Y)|_{U}) \ar[r]\ar@{^{(}->}[u]\ar@{^{(}->}[ru]& K_{0}^{top}(G(Y)|_{\partial\beta Y}) \ar[r]\ar[u]^{\mu_{bdry}}& K_{0}^{top}(U^{c}\rtimes \Gamma) \ar[r]\ar[u]^{\ucong}& K_{1}^{top}(X \times X)\ar[u]
}
\end{equation*}

\begin{remark}
We justify the diagonal inclusion of $K_{0}^{top}(G(Y)|_{U})$ into $K_{0}(C^{*}_{r}(G(Y)|_{\partial\beta Y})$. This follows as the groupoid $G(Y)|_{U}$ is amenable, and hence the assembly map is an isomorphism.   The algebra $C^{*}_{r}(G(Y)_{U})= \bigoplus_{i}M_{\vert X_{i}\vert} \otimes C(\partial\beta Y)$ injects into the product $\prod_{i} M_{\vert X_{i}\vert} \otimes C(\partial\beta Y)$ at the level of K-theory and this inclusion factors through in the inclusion into the kernel of $\pi$ and into $C^{*}_{r}(G(Y)|_{\partial\beta Y})$. These maps provide enough information to conclude injectivity of the assembly map $\mu_{bdry}$.
\end{remark}

A diagram chase under the assumption that the map $\mu_{bdry}$ is surjective quickly yields a contradiction, whence we have:

\begin{corollary}
The assembly map $\mu_{bdry}$ associated to $Y=\sqcup_{i,j}Cay(\frac{\Gamma}{N_{i}})$ is not surjective but is injective.\qed
\end{corollary}

\bibliography{ref.bib}

\begin{thebibliography}{CTWY08}

\bibitem[BR84]{MR745358}
Jean-Camille Birget and John Rhodes.
\newblock Almost finite expansions of arbitrary semigroups.
\newblock {\em J. Pure Appl. Algebra}, 32(3):239--287, 1984.

\bibitem[CTWY08]{MR2419930}
Xiaoman Chen, Romain Tessera, Xianjin Wang, and Guoliang Yu.
\newblock Metric sparsification and operator norm localization.
\newblock {\em Adv. Math.}, 218(5):1496--1511, 2008.

\bibitem[Exe08]{MR2419901}
Ruy Exel.
\newblock Inverse semigroups and combinatorial {$C\sp \ast$}-algebras.
\newblock {\em Bull. Braz. Math. Soc. (N.S.)}, 39(2):191--313, 2008.

\bibitem[FW95]{MR1388300}
Steven~C. Ferry and Shmuel Weinberger.
\newblock A coarse approach to the {N}ovikov conjecture.
\newblock In {\em Novikov conjectures, index theorems and rigidity, {V}ol.\ 1
  ({O}berwolfach, 1993)}, volume 226 of {\em London Math. Soc. Lecture Note
  Ser.}, pages 147--163. Cambridge Univ. Press, Cambridge, 1995.

\bibitem[GTY11]{MR2764895}
Erik Guentner, Romain Tessera, and Guoliang Yu.
\newblock Operator norm localization for linear groups and its application to
  {$K$}-theory.
\newblock {\em Adv. Math.}, 226(4):3495--3510, 2011.

\bibitem[GWY08]{MR2431253}
Guihua Gong, Qin Wang, and Guoliang Yu.
\newblock Geometrization of the strong {N}ovikov conjecture for residually
  finite groups.
\newblock {\em J. Reine Angew. Math.}, 621:159--189, 2008.

\bibitem[Hig99]{higsonpreprint}
Nigel Higson.
\newblock Counterexamples to the coarse {B}aum-{C}onnes conjecture.
\newblock {\em Preprint}, 1999.

\bibitem[Hig00]{MR1779613}
N.~Higson.
\newblock Bivariant {$K$}-theory and the {N}ovikov conjecture.
\newblock {\em Geom. Funct. Anal.}, 10(3):563--581, 2000.

\bibitem[HLS02]{MR1911663}
N.~Higson, V.~Lafforgue, and G.~Skandalis.
\newblock Counterexamples to the {B}aum-{C}onnes conjecture.
\newblock {\em Geom. Funct. Anal.}, 12(2):330--354, 2002.

\bibitem[HR95]{MR1388312}
Nigel Higson and John Roe.
\newblock On the coarse {B}aum-{C}onnes conjecture.
\newblock In {\em Novikov conjectures, index theorems and rigidity, {V}ol.\ 2
  ({O}berwolfach, 1993)}, volume 227 of {\em London Math. Soc. Lecture Note
  Ser.}, pages 227--254. Cambridge Univ. Press, Cambridge, 1995.

\bibitem[HR00]{MR1817560}
Nigel Higson and John Roe.
\newblock {\em Analytic {$K$}-homology}.
\newblock Oxford Mathematical Monographs. Oxford University Press, Oxford,
  2000.
\newblock Oxford Science Publications.

\bibitem[KL04]{MR2041539}
J.~Kellendonk and Mark~V. Lawson.
\newblock Partial actions of groups.
\newblock {\em Internat. J. Algebra Comput.}, 14(1):87--114, 2004.

\bibitem[K{\"o}n90]{MR1035708}
D{\'e}nes K{\"o}nig.
\newblock {\em Theory of finite and infinite graphs}.
\newblock Birkh\"auser Boston Inc., Boston, MA, 1990.
\newblock Translated from the German by Richard McCoart, With a commentary by
  W. T. Tutte and a biographical sketch by T. Gallai.

\bibitem[KS02]{MR1900993}
Mahmood Khoshkam and Georges Skandalis.
\newblock Regular representation of groupoid {$C^*$}-algebras and applications
  to inverse semigroups.
\newblock {\em J. Reine Angew. Math.}, 546:47--72, 2002.

\bibitem[Law10]{MR2672179}
M.~V. Lawson.
\newblock Compactable semilattices.
\newblock {\em Semigroup Forum}, 81(1):187--199, 2010.

\bibitem[LMS06]{MR2221438}
Mark~V. Lawson, Stuart~W. Margolis, and Benjamin Steinberg.
\newblock Expansions of inverse semigroups.
\newblock {\em J. Aust. Math. Soc.}, 80(2):205--228, 2006.

\bibitem[OOY09]{MR2568691}
Herv{\'e} Oyono-Oyono and Guoliang Yu.
\newblock {$K$}-theory for the maximal {R}oe algebra of certain expanders.
\newblock {\em J. Funct. Anal.}, 257(10):3239--3292, 2009.

\bibitem[Pat99]{MR1724106}
Alan L.~T. Paterson.
\newblock {\em Groupoids, inverse semigroups, and their operator algebras},
  volume 170 of {\em Progress in Mathematics}.
\newblock Birkh\"auser Boston Inc., Boston, MA, 1999.

\bibitem[Roe03]{MR2007488}
John Roe.
\newblock {\em Lectures on coarse geometry}, volume~31 of {\em University
  Lecture Series}.
\newblock American Mathematical Society, Providence, RI, 2003.

\bibitem[STY02]{MR1905840}
G.~Skandalis, J.~L. Tu, and G.~Yu.
\newblock The coarse {B}aum-{C}onnes conjecture and groupoids.
\newblock {\em Topology}, 41(4):807--834, 2002.

\bibitem[Tu99]{MR1703305}
Jean-Louis Tu.
\newblock La conjecture de {B}aum-{C}onnes pour les feuilletages moyennables.
\newblock {\em $K$-Theory}, 17(3):215--264, 1999.

\bibitem[Tu00]{MR1798599}
Jean-Louis Tu.
\newblock The {B}aum-{C}onnes conjecture for groupoids.
\newblock In {\em {$C^*$}-algebras ({M}\"unster, 1999)}, pages 227--242.
  Springer, Berlin, 2000.

\bibitem[Wan07]{MR2363697}
Qin Wang.
\newblock Remarks on ghost projections and ideals in the {R}oe algebras of
  expander sequences.
\newblock {\em Arch. Math. (Basel)}, 89(5):459--465, 2007.

\bibitem[WY12]{explg1}
Rufus Willett and Guoliang Yu.
\newblock Higher index theory for certain expanders and {G}romov monster
  groups, {I}.
\newblock {\em Adv. Math.}, 229(3):1380--1416, 2012.

\end{thebibliography}

\end{document}